\documentclass{article}
\title{Triviality Theorems for Yetter-Drinfel'd Hopf Algebras}
\author{Yorck Sommerh\"auser}
\date{}

\usepackage{amsmath}
\usepackage{theorem}
\usepackage{amssymb}
\usepackage{latexsym}

\makeatletter
\renewcommand{\subsection}{\@startsection{subsection}{2}{0em}%
{\baselineskip}{-0em}{\bfseries\normalsize}}
\makeatother

\makeatletter
\newcommand{\listofdefinitions}{\@starttoc{def}}

\newcommand{\l@definition}[2]{\par\noindent#1 {\itshape #2}}
\makeatother

\newcounter{num}
\newenvironment{pflist}{\begin{list}{(\arabic{num})}{\usecounter{num} \leftmargin0cm \itemindent5pt}}{\end{list}}

\theoremstyle{plain}
\theorembodyfont{\rmfamily}

\newtheorem{thm}{Theorem}

\newtheorem{prop}[thm]{Proposition}
\newtheorem{lemma}[thm]{Lemma}
\newtheorem{corollary}[thm]{Corollary}
\newtheorem{pf}{Proof.}

\newtheorem{defn}[thm]{Definition}

\theoremstyle{break}
\theorembodyfont{\rmfamily}
\newtheorem{propb}[thm]{Proposition}
\newtheorem{defb}[thm]{Definition}
\newtheorem{lemb}[thm]{Lemma}
\newtheorem{thmb}[thm]{Theorem}
\newtheorem{corb}[thm]{Corollary}

\newcommand{\qed}{$\Box$}

\newcommand{\End}{\operatorname{End}}
\newcommand{\Aut}{\operatorname{Aut}}
\newcommand{\Span}{\operatorname{Span}}
\newcommand{\Hom}{\operatorname{Hom}}
\newcommand{\id}{\operatorname{id}}

\newcommand{\op}{\scriptstyle \operatorname{op}}
\newcommand{\cop}{\scriptstyle \operatorname{cop}}

\def\1{{(1)}}
\def\2{{(2)}}
\def\3{{(3)}}
\def\4{{(4)}}
\def\5{{(5)}}

\def\o{\otimes}
\def\ra{\rightarrow}

\def\da{\Delta_{A}}

\def\dh{\Delta_{H}}
\def\ea{\varepsilon_{A}}
\def\eb{\varepsilon_{B}}
\def\eh{\varepsilon_{H}}

\def\sa{S_{A}}

\def\sh{S_{H}}

\def\ioa{\iota_A}
\def\lma{\lambda_A}

\def\Lma{{\Lambda_A}{}}

\def\lmb{\lambda_B}

\def\Lmb{\Lambda_B}

\def\lmh{\lambda_H}

\def\Lmh{\Lambda_H}

\def\H{1_H}
\def\A{1_A}
\def\B{1_B}

\def\N{{\mathbb N}}
\def\E{{E}}

\begin{document}

\maketitle

\begin{abstract}
\hspace{-5mm}Under suitable assumptions on the base field, we prove that a commutative semisimple Yetter-Drinfel'd Hopf algebra over a finite abelian group is trivial, i.e., is an ordinary Hopf algebra, if its dimension is relatively prime to the order of the finite abelian group. Furthermore, we prove that a finite-dimensional cocommutative cosemisimple Yetter-Drinfel'd Hopf algebra contains a trivial Yetter-Drinfel'd Hopf subalgebra of dimension greater than one, at least if the Yetter-Drinfel'd Hopf algebra itself has dimension greater than one. 
\end{abstract}


\section*{Introduction} \label{Sec:Introd}
\addcontentsline{toc}{section}{Introduction}
Usually, a Yetter-Drinfel'd Hopf algebra is not an ordinary Hopf algebra. The difference between the two notions is caused by the fact that the coproduct of a Yetter-Drinfel'd Hopf algebra is not an algebra homomorphism with respect to the canonical algebra structure on the second tensor power of the Yetter-Drinfel'd Hopf algebra. Rather, it is an algebra homomorphism with respect to a modified algebra structure on the second tensor power constructed via a special quasisymmetry that is characteristic for Yetter-Drinfel'd modules. However, it may happen that, for the specific Yetter-Drinfel'd Hopf algebra under consideration, this quasisymmetry coincides on the second tensor power of the algebra with the ordinary flip of tensor factors. Then the algebra structure is in fact not modified, and the Yetter-Drinfel'd Hopf algebra is an ordinary Hopf algebra.

In this case, the Yetter-Drinfel'd Hopf algebra is called trivial (cf.~\cite{SoYp}, Def.~1.1, p.~8). According to a result of P. Schauenburg, the behavior of the quasisymmetry just described characterizes this property: If a Yetter-Drinfel'd Hopf algebra is an ordinary Hopf algebra, then the quasisymmetry for Yetter-Drinfel'd modules coincides on the second tensor power of the algebra with the ordinary flip of tensor factors (cf.~\cite{Schau}, Cor.~2, p.~262; see also~\cite{SoYp}, Prop.~1.1, p.~8).

The first main result of this article is that, in a certain situation, every Yetter-Drinfel'd Hopf algebra is trivial. In this situation, we consider a semisimple commutative Yetter-Drinfel'd Hopf algebra~$A$ over the group ring~$K[G]$ of a finite abelian group~$G$, where~$K$ is a field of characteristic zero, and prove in Paragraph~\ref{TrivThm} the following triviality theorem:
\begin{thm}
If~$\dim(A)$ and~$|G|$ are relatively prime, then $A$ is trivial.
\end{thm}
This result was known in the case where~$|G|$ is prime (cf.~\cite{SoYp}, Cor.~6.7, p.~100).

The second main result of this article is concerned with a finite-dimensional cocommutative cosemisimple Yetter-Drinfel'd Hopf algebra~$A$ over the group ring~$K[G]$ of a finite abelian group~$G$, where~$K$ is an algebraically closed field of characteristic zero. We also prove in Paragraph~\ref{TrivThm} that, although~$A$ itself need not be trivial, it contains at least a trivial part:
\begin{prop}
If $\dim(A)>1$, then $A$ contains a trivial Yetter-Drinfel'd Hopf subalgebra~$B$ with $\dim(B)>1$.
\end{prop}

We develop the rather involved general theory of such Yetter-Drinfel'd Hopf algebras here only to the extent that is necessary to prove these two results. Section~\ref{Sec:Prelim} contains a brief, but nonetheless very important, summary of basic facts about Yetter-Drinfel'd Hopf algebras. More detailed treatments of this material can be found in \cite{SoDefUnivEinh}, \cite{SoDevEnvAlg}, \cite{SoRib}, \cite{SoDiss}, or~\cite{SoYp}. Section~\ref{Sec:Prelim} ends with Theorem~\ref{Int}, a version of the Nichols-Zoeller theorem for two Yetter-Drinfel'd Hopf algebras that are not necessarily defined over the same Hopf algebra. In Section~\ref{Sec:FinAbelGroups}, we discuss Yetter-Drinfel'd Hopf algebras over group rings of finite abelian groups. In particular, we introduce there a certain alternative description of the coaction that will be used in the entire article. This brief section also contains an important discussion of a way to modify the finite abelian group over which the Yetter-Drinfel'd Hopf algebra is defined. The main part of the article is Section~\ref{Sec:CommCase}, which deals with commutative semisimple Yetter-Drinfel'd Hopf algebras over group rings of finite abelian groups. As we will indicate during the discussion, most of the material presented there generalizes facts that were established in Chapter~6 of~\cite{SoYp} for the case in which the finite abelian group has prime order. This final section of the article ends with the proof of the two main results stated above.

While the precise assumptions that are used throughout each section are listed in its first paragraph, there are some conventions that are common to all sections. The base field is always denoted by~$K$, and while it is arbitrary in Section~\ref{Sec:Prelim}, we assume that it is algebraically closed of characteristic zero in Section~\ref{Sec:FinAbelGroups} and Section~\ref{Sec:CommCase} unless this is explicitly stated otherwise, which happens only in a minor comment in Paragraph~\ref{TrivThm}. The multiplicative group of invertible elements in the field~$K$ is denoted by~$K^\times := K\setminus\{0\}$. All vector spaces are defined over~$K$, and all unadorned tensor products are taken over~$K$. The dual of a vector space~$V$ is denoted by $V^*:=\Hom_K(V,K)$, and the transpose of a linear map~$f$, i.e., the induced map between the dual spaces, is denoted by~$f^*$. The symbol $\perp$ will be used in four different meanings, two for vector spaces and two for groups. In this article, a character is a one-dimensional character, i.e., a group homomorphism to the multiplicative group~$K^\times$ in the case of a group character, or an algebra homomorphism to the base field~$K$ in the case of the character of an algebra. All rings, and therefore especially all algebras, are assumed to have a unit element, and ring and algebra homomorphisms are assumed to preserve unit elements. Unless stated otherwise, a module is a left module. The cardinality of a set~$X$ is denoted by~$|X|$. The symbol $\subset$ denotes non-strict inclusion, so that we have $X \subset X$ for every set~$X$. Also, we use the so-called Kronecker symbol~$\delta_{ij}$, which is equal to~$1$ if $i=j$ and equal to~$0$ otherwise. With respect to enumeration, we use the convention that propositions, definitions, and similar items are referenced by the paragraph in which they occur; i.e., a reference to Proposition~1.1 refers to the unique proposition in Paragraph~1.1. 

The first main result of this article was presented at the AMS Fall Eastern Sectional Meeting in Halifax in October~2014, while the second main result was presented at the Joint Mathematics Meeting in San Antonio in January~2015. The author thanks the organizers of these conferences for the invitation. He also thanks the Department of Mathematics at SUNY Buffalo for a visiting appointment during which most of this article was written.

\section{Preliminaries} \label{Sec:Prelim}
\subsection[Left Yetter-Drinfel'd modules]{} \label{LeftYetMod}
Although we will limit our considerations to special situations soon, let us briefly recall the case of a general Hopf algebra~$H$ with coproduct~$\dh$, counit~$\eh$, and antipode~$\sh$. A left-left Yetter-Drinfel'd module, or left Yetter-Drinfel'd module for short, is a left $H$-module~$V$ that is simultaneously a left $H$-comodule in such a way that these two structures are compatible in the sense that
$$\delta(h.v) = h_{\1} v^{\1} \sh(h_{\3}) \o h_{\2}.v^{\2}$$
for all $h \in H$ and all $v \in V$. Here we have used Heyneman-Sweedler sigma notation in the form 
$$\dh(h) = h_{\1} \o h_{\2} \in H \o H \qquad \qquad
\delta(v) = v^{\1} \o v^{\2} \in H \o V$$ 
for the coproduct~$\dh$ and the coaction~$\delta$, respectively. The dot denotes the module action. Although right Yetter-Drinfel'd modules will be defined below, we always mean left Yetter-Drinfel'd modules if no side is specified. We note already at this point that, in the case where~$H$ is commutative and cocommutative, which will be the case in Section~\ref{Sec:FinAbelGroups} and Section~\ref{Sec:CommCase}, the Yetter-Drinfel'd condition stated above reduces to
$$\delta(h.v) = v^{\1} \o h.v^{\2}$$
which coincides with the compatibility condition for left-left dimodules (cf.~\cite{Long}, Def.~3.1, p.~575).

The base field~$K$ is a Yetter-Drinfel'd module if endowed with the trivial module structure
$h.\xi := \eh(h) \xi$ and the trivial comodule structure $\delta(\xi) := \H \o \xi$, for $\xi \in K$. Moreover, the tensor product of two Yetter-Drinfel'd modules, together with the usual diagonal module structure and the codiagonal comodule structure, is again a Yetter-Drinfel'd module. If~$V$ and~$W$ are two Yetter-Drinfel'd modules, the two possible ways of forming their tensor product are related via the quasisymmetry
$$\sigma_{V,W}: V \o W \rightarrow W \o V,~v \o w \mapsto (v^{\1}.w \o v^{\2})$$
This map is both $H$-linear and colinear. If the antipode of~$H$ is bijective, the quasisymmetry is also bijective, and its inverse is then given by the formula 
$\sigma_{V,W}^{-1}(w \o v) = v^{\2} \o \sh^{-1}(v^{\1}).w$.

\subsection[Right Yetter-Drinfel'd modules]{} \label{RightYetMod}
Similarly, a right-right Yetter-Drinfel'd module, or right Yetter-Drinfel'd module for short, is a right $H$-module~$V$ that is simultaneously a right \mbox{$H$-comod}\-ule in such a way that these two structures are compatible in the sense that
$$\delta(v.h) = v^{\1}.h_{\2} \o \sh(h_{\1}) v^{\2} h_{\3}$$
for all $h \in H$ and all $v \in V$. It is important to note that here, in contrast to the previous case, $v.h \in V$ and $\delta(v) = v^{\1} \o v^{\2} \in V \o H$. As right $H$-modules are just left modules over the opposite algebra~$H^{\op}$ and right $H$-comodules are just left comodules over the coopposite coalgebra~$H^{\cop}$, right Yetter-Drinfel'd modules over~$H$ are just left Yetter-Drinfel'd modules over~$H^{\op\cop}$. It is also possible to define left-right and right-left Yetter-Drinfel'd modules, but they will not be used in the sequel. As before, in the case where~$H$ is commutative and cocommutative, the right Yetter-Drinfel'd condition coincides with the right-right dimodule condition
$$\delta(v.h) = v^{\1}.h \o  v^{\2}$$
Again, the base field~$K$ is a right Yetter-Drinfel'd module if endowed with the trivial module structure
$\xi.h := \eh(h) \xi$ and the trivial comodule structure $\delta(\xi) := \xi \o \H$, for $\xi \in K$, and the tensor product of two right Yetter-Drinfel'd modules~$V$ and~$W$, together with the usual diagonal module structure and the codiagonal comodule structure, is again a right Yetter-Drinfel'd module. The formula for the quasisymmetry is $\sigma_{V,W}(v \o w) = w^{\1} \o v.w^{\2}$ in this case. As before, the quasisymmetry is bijective with inverse 
$$\sigma_{V,W}^{-1}(w \o v) = v.\sh^{-1}(w^{\2}) \o w^{\1}$$ 
if the antipode~$\sh$ is bijective. It should be noted that it is necessary to reverse the tensor factors when considering right Yetter-Drinfel'd modules as left Yetter-Drinfel'd modules over~$H^{\op\cop}$: If the right Yetter-Drinfel'd modules~$V$ and~$W$ are viewed as left Yetter-Drinfel'd modules over~$H^{\op\cop}$, their tensor \mbox{product $V \o W$} over~$H^{\op\cop}$ does not correspond to their tensor product~$V \o W$ over~$H$, but rather to the tensor product~$W \o V$ over~$H$, viewed as a left Yetter-Drinfel'd module over~$H^{\op\cop}$.  This reversion of the tensor factors is also compatible with the quasisymmetries. In the language of category theory, this means that we have constructed a (non-strict) quasisymmetric monoidal functor from the category of right Yetter-Drinfel'd modules over~$H$ to the reverse category of the category of left Yetter-Drinfel'd modules over~$H^{\op\cop}$ (cf.~\cite{JoyStreet}, Ex.~2.5, p.~39; \cite{Kas}, Exerc.~XIII.6.4, p.~337).

In our treatment, right Yetter-Drinfel'd modules arise as duals of left Yetter-Drinfel'd modules as a consequence of the following proposition:
\begin{prop}
If $V$ is a finite-dimensional left Yetter-Drinfel'd module, then the dual space~$V^*$ is in a unique way a right Yetter-Drinfel'd module such that the natural
pairing
$$\langle \cdot , \cdot \rangle : V \times V^* \rightarrow K, 
(v,\varphi) \mapsto \langle v, \varphi \rangle := \varphi(v)$$
is a Yetter-Drinfel'd form in the sense that we have
\begin{enumerate}
\item $\langle h.v, \varphi \rangle = \langle v, \varphi.h \rangle$
\item $\langle v,\varphi^{\1}\rangle \varphi^{\2} = v^{\1} \langle v^{\2},\varphi \rangle$
\end{enumerate}
\end{prop}

If $W$ is another finite-dimensional left Yetter-Drinfel'd module, then the pairing between 
$V \o W$ and $V^* \o W^*$ defined by 
$$\langle v \otimes w, \varphi \o \varphi' \rangle = 
\langle v, \varphi \rangle \langle w, \varphi' \rangle$$
is also a Yetter-Drinfel'd form, and the quasisymmetries are mutually adjoint with respect to this pairing; i.e., we have
$$\langle \sigma_{V, W}(v \otimes w), \varphi' \o \varphi \rangle = 
\langle v \otimes w, \sigma_{W^*, V^*}(\varphi' \o \varphi) \rangle$$
for $v \in V$, $w \in W$, $\varphi \in V^*$, and  $\varphi' \in W^*$. In the language of category theory, this means that we have constructed a (non-strict) quasisymmetric monoidal functor from the category of finite-dimensional left Yetter-Drinfel'd modules to the category of finite-dimensional right Yetter-Drinfel'd modules. More precisely, since this functor is contravariant, it is a functor to the opposite category of the category of finite-dimensional right Yetter-Drinfel'd modules (cf.~\cite{Tur}, Thm.~II.5.4, p.~112). Proofs of these statements, which are not difficult, are briefly indicated in \cite{SoDevEnvAlg}, Subsec.~2.4 and Subsec.~2.5, p.~37; more detailed proofs can be found in the corresponding parts of~\cite{SoDefUnivEinh}.

\subsection[Yetter-Drinfel'd Hopf algebras]{} \label{YetDrinfHopf}
A left (or right) Yetter-Drinfel'd algebra is a left (or right) Yetter-Drinfel'd module which is simultaneously an associative algebra with unit in such a way that the multiplication map
$$\mu_A: A \o A \rightarrow A,~a \o a' \mapsto aa'$$
and the unit map $K \rightarrow A,~\xi \mapsto \xi \A$ are both $H$-linear and colinear. If~$A'$ is a second left (or right) Yetter-Drinfel'd algebra, then the tensor product~$A \o A'$ becomes an associative algebra with unit with respect to the multiplication
$$\mu_{A \o A'} := (\mu_{A} \o \mu_{A'}) \circ (\id_A \o \sigma_{A' \o A} \o \id_{A'})$$
The difference to the canonical tensor algebra structure consists in the fact that this algebra structure uses the quasisymmetry~$\sigma_{A' \o A}$ instead of the usual flip of tensor factors. We denote this tensor product by~$A \hat{\o} A'$ if we want to emphasize that it carries this algebra structure in contrast to the canonical tensor algebra structure. For left Yetter-Drinfel'd algebras, the explicit form of this product is
$$(a \o a')(b \o b') = a(a'{}^\1.b) \o a'{}^\2 b'$$
which shows that this multiplication is a special case of a left smash product of the module algebra~$A$ and the comodule algebra~$A'$ (cf.~\cite{Beattie}, Def.~1.2, p.~24; see also \cite{Doi1}, Rem.~(1.3), p.~374 and \cite{SoYp}, Par.~1.7, p.~15). 

Analogously, a left (or right) Yetter-Drinfel'd coalgebra is a left (or right) Yetter-Drinfel'd module which is simultaneously a coalgebra in such a way that the comultiplication $\da: A  \rightarrow A \o A$ and the counit $\ea: A \rightarrow K$ are both \mbox{$H$-linear} and colinear. A left (or right) Yetter-Drinfel'd bialgebra is a left (or right) Yetter-Drinfel'd algebra which is simultaneously a Yetter-Drinfel'd coalgebra with the property that both the coproduct~$\da$ and the counit~$\ea$ are algebra homomorphisms. Here, when we say that~$\da$ is an algebra homomorphism, we refer to the algebra structure~$A \hat{\o} A$ described above, and not to the canonical tensor algebra structure; i.e., we view the coproduct as a map $\da: A  \rightarrow A \hat{\o} A$. Finally, a left (or right) Yetter-Drinfel'd Hopf algebra is a left (or right) Yetter-Drinfel'd bialgebra together with an $H$-linear and colinear map~$\sa: A \rightarrow A$, called the antipode, that satisfies
$$\sa(a_\1) a_\2 = a_\1 \sa(a_\2) = \ea(a) \A$$
for all $a \in A$. As in the case of ordinary Hopf algebras, such an antipode satisfies $\sa(\A) = \A$ and 
$\ea \circ \sa = \ea$, but in the compatibility with product and coproduct the flip of tensor factors must again be replaced by the quasisymmetry, namely in the form
$$\sa \circ \mu_A = \mu_A \circ (\sa \o \sa) \circ \sigma_{A,A} \qquad \qquad
\da \circ \sa = \sigma_{A,A} \circ (\sa \o \sa) \circ \da$$
(cf.~\cite{Ly1}, Prop.~3.1, p.~270; \cite{M}, \S~10.5, Eq.~(10.5.9), p.~205; \cite{SoRib}, Par.~2.5, p.~427).
As for Hopf algebras, the antipode of a finite-dimensional Yetter-Drinfel'd Hopf algebra is bijective (cf.~\cite{Doi2}, Thm.~3, p.~3066; \cite{SoRib}, Prop.~2.10, p.~432).

It is obvious from our definitions that a Yetter-Drinfel'd Hopf algebra is in general not an ordinary Hopf algebra, because the algebra structure on~$A \hat{\o} A$ is in general different from the canonical tensor algebra structure. However, it is also obvious from our definitions that these two algebra structures coincide if it happens that $\sigma_{A \o A}(a \o a') = a' \o a$ for all $a,a' \in A$, in which case~$A$ is indeed an ordinary Hopf algebra. As observed by P.~Schauenburg, the converse of this statement is also true: If a Yetter-Drinfel'd Hopf algebra is an ordinary Hopf algebra, then we have
$\sigma_{A \o A}(a \o a') = a' \o a$ for all $a,a' \in A$ (cf.~\cite{Schau}, Cor.~2, p.~262; see also~\cite{SoYp}, Prop.~1.1, p.~8). We call such Yetter-Drinfel'd Hopf algebras trivial:
\begin{defn}
A left (or right) Yetter-Drinfel'd Hopf algebra~$A$ is called trivial if we have $\sigma_{A,A}(a \o a') = a' \o a$ for all $a,a' \in A$.
\end{defn}

One of the main results of this article, namely Theorem~\ref{TrivThm}, describes a situation that forces a Yetter-Drinfel'd Hopf algebra to be trivial.

We have already discussed in Paragraph~\ref{RightYetMod} how the dual space of a finite-dimensional left Yetter-Drinfel'd module becomes a right Yetter-Drinfel'd module. If the left Yetter-Drinfel'd module was a Yetter-Drinfel'd Hopf algebra, the arising right Yetter-Drinfel'd module is also a Yetter-Drinfel'd Hopf algebra:
\begin{prop}
If $A$ is a finite-dimensional left Yetter-Drinfel'd Hopf algebra, then the dual space $A^*$ is in a unique way a right Yetter-Drinfel'd Hopf algebra such that the natural pairing described above is not only a Yetter-Drinfel'd form, but also a bialgebra form, which means that we have
\begin{enumerate}
\item $\langle a \otimes a', \Delta_{A^*}(\varphi)\rangle = \langle aa', \varphi \rangle$

\item $\langle a, \varphi \varphi'\rangle = \langle \Delta_{A}(a), \varphi \otimes \varphi'\rangle$

\item $\langle 1_A, \varphi\rangle = \varepsilon_{A^*}(\varphi),\, \langle a, 1_{A^*} \rangle = \ea(a)$
\end{enumerate}
where the pairing of the second tensor powers is defined as in Paragraph~\ref{RightYetMod}.
\end{prop}
This proposition appears in several places, among them \cite{SoDevEnvAlg}, Subsec.~2.6, p.~37; the corresponding part of~\cite{SoDefUnivEinh} contains further details. More material about Yetter-Drinfel'd Hopf algebras can be found in \cite{M}, \cite{RadfHopf}, \cite{SoRib}, or \cite{SoYp}, and some of the results established there will be used in the sequel.

We have also discussed in Paragraph~\ref{RightYetMod} that right Yetter-Drinfel'd modules over~$H$ are just left Yetter-Drinfel'd modules over $H^{\op \cop}$, up to a reversion of the order of tensor products. For Yetter-Drinfel'd Hopf algebras, this implies the following:
\enlargethispage{3mm}
\begin{lemb} 
\begin{enumerate}
\item If $A$ is a left Yetter-Drinfel'd Hopf algebra over $H$, then the opposite algebra and coopposite coalgebra~$A^{\op \cop}$ is a right Yetter-Drinfel'd Hopf algebra over $H^{\op \cop}$.

\item If $A$ is a right Yetter-Drinfel'd Hopf algebra over $H$, then $A^{\op \cop}$ is a left Yetter-Drinfel'd Hopf algebra over $H^{\op \cop}$.
\end{enumerate}
\end{lemb}

This result is already stated in \cite{SoDevEnvAlg}, Lem.~2.4, p.~39; again, the corresponding part of \cite{SoDefUnivEinh} contains a slightly more detailed treatment. The result is also stated in \cite{SoRib}, Lem.~2.5, p.~427 and \cite{SoYp}, Lem.~1.2, p.~9.

\subsection[The dual Hopf algebra]{} \label{DualHopf}
In Proposition~\ref{RightYetMod} and Proposition~\ref{YetDrinfHopf}, we have dualized Yetter-Drinfel'd modules and Yetter-Drinfel'd Hopf algebras, respectively. But if the underlying Hopf algebra~$H$ is finite-dimensional, the process of dualization can also be applied to~$H$: If $H$ is finite-dimensional and $V$ is a left Yetter-Drinfel'd module over $H$, we know from \cite{M}, Lem.~1.6.4, p.~11 that the left comodule structure determines a right $H^*$-module structure, namely the module structure
$$V \o H^* \ra V, \: v \o \varphi \mapsto \varphi(v^\1) v^\2$$
Dually, the left $H$-module structure determines a right $H^*$-comodule structure
$$\delta_*: V \ra V \o H^*, \: v \mapsto v^{[1]} \o v^{[2]}$$
such that the analogous condition $(v^{[2]}(h)) v^{[1]} = h.v$
is satisfied for all $v \in V$ and $h \in H$, where we have used square brackets to distinguish the $H^*$-coaction from the $H$-coaction. It is then straightforward to verify that $V$ is a right Yetter-Drinfel'd module over~$H^*$.

By construction, a mapping between two Yetter-Drinfel'd modules is linear and colinear with respect to~$H$ if and only if it is linear and colinear with respect to~$H^*$. This dualization process commutes with taking the tensor product of two Yetter-Drinfel'd modules, and the quasisymmetry on the tensor product is the same before and after the dualization, because the preceding compatibility conditions imply that
$$w^{[1]} \o v.w^{[2]} = w^{[1]} \o w^{[2]}(v^{\1}) v^{\2} = v^{\1}.w \o v^{\2}$$ 
for elements of two left Yetter-Drinfel'd modules~$V$ and~$W$. In the language of category theory, this means that we have constructed a strict quasisymmetric monoidal functor from the category of left Yetter-Drinfel'd modules over~$H$ to the category of right Yetter-Drinfel'd modules over~$H^*$. Therefore, if $A$ is a left Yetter-Drinfel'd Hopf algebra over~$H$, it can also be considered as a right Yetter-Drinfel'd Hopf algebra over~$H^*$. Similarly, if $A$ is a right Yetter-Drinfel'd Hopf algebra over~$H$, it can also be considered as a left Yetter-Drinfel'd Hopf algebra over~$H^*$, as can be deduced from Lemma~\ref{YetDrinfHopf}. This topic is also discussed in \cite{SoRib}, Lem.~2.4, p.~427 and \cite{SoYp}, Lem.~1.3, p.~10.

\subsection[Integrals]{} \label{Int}
We will need some results from the theory of integrals for Yetter-Drinfel'd Hopf algebras. The basic facts about integrals for Yetter-Drinfel'd Hopf algebras were already obtained by D.~E.~Radford in \cite{RadfProj}, although some results there were proved only under the additional assumption that the Yetter-Drinfel'd Hopf algebra is an ordinary Hopf algebra; i.e., that it is trivial in the sense of Definition~\ref{YetDrinfHopf}. The theory was afterwards developed further by several authors, for example in the articles \cite{Doi2}, \cite{FMS}, \cite{Ly1}, \cite{Ly2}, and \cite{SoRib}. Here, we only briefly collect those results that will become relevant in the sequel, following the treatment in~\cite{SoRib}:
\begin{prop}
Suppose that~$A$ is a semisimple left Yetter-Drinfel'd Hopf algebra over the semisimple Hopf algebra~$H$, and  that the base field~$K$ has characteristic zero. Then the following holds:
\begin{enumerate}
\item 
There is a unique two-sided integral~$\Lma \in A$ satisfying~$\ea(\Lma)=1$.

\item 
There is a unique two-sided integral~$\lma \in A^*$ satisfying~$\lma(\Lma)=1$.

\item 
$A$ is finite-dimensional. 

\item  
$\lma$ is the character of the regular representation. 

\item 
$\lma(\A)=\dim(A)$

\item 
$\sa(\Lma) = \Lma$ and $\lma \circ \sa = \lma$

\item 
$\lma$ is a Frobenius homomorphism with corresponding Casimir element 
$$\Lma_\1 \o \sa(\Lma_\2) = \sa(\Lma_\1) \o \Lma_\2$$
i.e., we have
$a = \lma(a \Lma_\1) \sa(\Lma_\2) = \sa(\Lma_\1) \lma(\Lma_\2 a)$
for all $a \in A$.

\item 
$\lma$ and~$\Lma$ are cocommutative.
\end{enumerate}
\end{prop}
\begin{pf}
\begin{pflist}
\item
The kernel of~$\ea$ is a two-sided ideal of codimension~$1$. Since~$A$ is semisimple, there is a two-sided ideal that complements this kernel; this complement is then one-dimensional. As in the case of Maschke's theorem for ordinary Hopf algebras (cf.~\cite{M}, Thm.~2.2.1, p.~20; \cite{RadfHopf}, Thm.~10.3.2, p.~298; \cite{Sw}, Thm.~5.1.8, p.~103), we can see that a nonzero element~$\Lma$ in this complement is a two-sided integral. Since~$\Lma$ is not contained in the kernel of~$\ea$, we have~$\ea(\Lma) \neq 0$; via rescaling, we can achieve that~$\ea(\Lma)=1$. By multiplying~$\Lma$ with a potentially different element with the same properties, we see that this element is unique, proving the first assertion.

\item
Since $H$ is semisimple, $H$ is finite-dimensional (cf.~\cite{RadfHopf}, Prop.~10.2.1, p.~293; \cite{Sw1}, Cor.~2.7, p.~330; \cite{Sw}, Chap.~V, Exerc.~4, p.~108). If $\Lmh \in H$ is a two-sided integral satisfying \mbox{$\eh(\Lmh)=1$}, then $\Lmb:=\Lma \o \Lmh$ is a right integral in the Radford biproduct~$B:=A \o H$ (cf.~\cite{RadfProj}, Prop.~3, p.~333). Since \mbox{$\ea(\Lma)\eh(\Lmh)=1$}, the Radford biproduct is semisimple, and therefore finite-dimensional by the result just cited. But this obviously means that~$A$ is finite-dimensional, which proves the third assertion.

\item
Now we can choose nonzero right integrals $\lma \in A^*$ and $\lmh \in H^*$ (cf.~\cite{Doi2}, Thm.~3, p.~3066; \cite{FMS}, Cor.~5.8, p.~4885; \cite{SoRib}, Prop.~2.10, p.~432). Then the formula
$$\lmb(a \o h) = \lma(a) \lmh(h)$$
defines a right integral~$\lmb$ on the Radford biproduct~$B=A \o H$ (cf.~\cite{RadfProj}, Prop.~4, p.~335). Since~$B$ is in this situation also cosemisimple (cf.~\cite{LR2}, Thm.~3.3, p.~276), we have
$$\lmb(\B) = \lma(\A) \lmh(\H) \neq 0$$
and therefore~$\lma(\A)  \neq 0$ and~$\lmh(\H) \neq 0$. By rescaling~$\lma$ and~$\lmh$, we can achieve that $\lma(\A) = \dim(A)$ and~$\lmh(\H) = \dim(H)$. Therefore also~$A$ is cosemisimple and both~$\lma$ and~$\lmh$ are two-sided integrals (cf.~\cite{SoRib}, Cor.~2.14, p.~438). Now it follows from \cite{LR1}, Thm.~4, p.~195 together with \cite{LR2}, Thm.~4.4, p.~279 that~$\lmb$ is the character of the (left or right) regular representation of~$B$. Since $(a \o \H)(a' \o h') = aa' \o h'$, it follows 
that~$\lmb(a \o \H)$ is $\dim(H)$ times the trace of the map $a' \mapsto aa'$. But since we clearly have
$\lmb(a \o \H) = \lma(a) \dim(H)$, this implies that this integral is the character of the regular representation. The fourth assertion will then follow if we can show that~$\lma$ also satisfies the normalization condition given in the second assertion, which we do next.

\item
From the first step of this proof, we know that~$\Lma$ is a centrally primitive idempotent that spans a one-dimensional two-sided ideal. It has therefore trace~$1$ in the regular representation. In view of the preceding step, this means that~$\lma(\Lma)=1$, proving the second assertion, as this normalization condition obviously determines the integral uniquely, and therefore also completing the proof of the fourth assertion.

\item
We have already obtained the fifth assertion in the proof of the fourth assertion above; it also follows immediately from the fourth assertion itself. In addition, it is a special case of a trace formula that appears in~\cite{AnSchn}, Thm.~7.3, p.~445 and~\cite{SoRib}, Par.~3.11, p.~447.

\item
By \cite{SoRib}, Prop.~2.12, p.~436, the antipode maps left integrals to right integrals and right integrals to left integrals. Since our integral is two-sided, $\sa(\Lma)$ is proportional to~$\Lma$. But as we stated in Paragraph~\ref{YetDrinfHopf}, we have $\ea \circ \sa = \ea$ and therefore $\ea(\sa(\Lma)) = \ea(\Lma) = 1$, which shows that the proportionality constant must be~$1$, so that $\sa(\Lma) = \Lma$. In view of Proposition~\ref{YetDrinfHopf} and Lemma~\ref{YetDrinfHopf}, the equation $\lma \circ \sa = \lma$ follows by applying this result to~$A^*$.

\item
The two forms of the Casimir element given in the seventh assertion follow from \cite{SoRib}, Prop.~2.10, p.~432 and Prop.~4.2, p.~449, respectively, because in our case the integrals are two-sided and the integral character~$\ioa$, the integral group element~$g_A$, and the right modular element~$a_A^R$ used there are \mbox{$\ioa=\eh$}, $g_A=\H$, and $a_A^R=\A$ (cf.~\cite{SoRib}, Prop.~2.14 and Cor.~2.14, p.~438; see also \cite{Doi2}, Thm.~3, p.~3066; \cite{FMS}, Thm.~5.6, p.~4883).

\item
For the eighth assertion, we note that the character of the regular representation is clearly cocommutative, so that~$\lma$ is cocommutative by the fourth assertion, i.e., satisfies \mbox{$\lma(ab) = \lma(ba)$} for all~$a,b \in A$. Alternatively, we know from \cite{SoRib}, Thm.~4.5, p.~454 that the negatively twisted Nakayama automorphism coincides with the square~$\sa^2$ of the antipode in our case, and then \cite{SoRib}, Prop.~4.6, p.~455 shows that the ordinary Nakayama automorphism is the composition of~$\sa^2$ with the ribbon transformation~$\theta_A$, defined in \cite{SoRib}, Par.~3.5, p.~441. But from the references cited in the third step, we know that the semisimplicity and cosemisimplicity of the Radford biproduct $B=A \o H$ implies that the square of its antipode is the identity, which entails by \cite{AnSchn}, Eq.~(4.4), p.~437 or \cite{SoRib}, Prop.~3.11, p.~447 that $\sa^2 \circ \theta_A$ is the identity. Thus the (ordinary) Nakayama automorphism is the identity, which means that~$\lma$ is cocommutative. By applying this result to~$A^*$, we get in view of Proposition~\ref{YetDrinfHopf} and Lemma~\ref{YetDrinfHopf} that~$\Lma$ is cocommutative, i.e., that $\Lma_\1 \o \Lma_\2 = \Lma_\2 \o \Lma_\1$.~\qed
\end{pflist}
\end{pf}

In the sequel, we will always use integrals that satisfy the two normalization conditions~$\ea(\Lma)=1$ and~$\lma(\Lma)=1$ appearing in this proposition. We note that the converse of its first assertion holds: If there exists a (left) integral~$\Lma$ satisfying~$\ea(\Lma)=1$, then~$A$ is semisimple (cf.~\cite{SoRib}, Prop.~2.14, p.~438). This shows that $A$ remains semisimple under extension of the base field; i.e., $A$ is separable (cf.~\cite{CR2}, Def.~(7.1), p.~142).

\begin{sloppypar}
The Nichols-Zoeller freeness theorem for ordinary Hopf algebras (cf.~\cite{M}, 
\mbox{Thm.~3.1.5}, p.~30; \cite{NZ}, Thm.~7, p.~384; \cite{RadfHopf}, Thm.~9.3.3, p.~281) has been generalized to Yetter-Drinfel'd Hopf algebras by B.~Scharfschwerdt (cf.~\cite{Scharf}, Thm.~2.2, p.~2483; see also \cite{ScharfDiss}; \cite{AnSchn}, Prop.~4.9, p.~438; \cite{SoYp}, Prop.~1.8, p.~16). We give here a slightly different version in which the smaller Yetter-Drinfel'd Hopf algebra is not necessarily defined over the same Hopf algebra as the larger one:
\end{sloppypar}
\begin{thm}
For a base field~$K$ of characteristic zero, we consider a semisimple left Yetter-Drinfel'd Hopf algebra~$A$ over the semisimple Hopf algebra~$H$ that contains a unital subalgebra~$B \subset A$ that is simultaneously a subcoalgebra. We assume that $B$ is also a left Yetter-Drinfel'd Hopf algebra over a potentially different semisimple Hopf algebra~$H'$. Then the following holds:
\begin{enumerate}
\item 
$B$ is semisimple.

\item 
$A$ is free as a left and right module over~$B$.

\item
$\dim(B)$ divides~$\dim(A)$.
\end{enumerate}
\end{thm}
\begin{pf}
\begin{pflist}
\item
In the third step of the proof of the preceding proposition, we have already described the integrals on the Radford biproduct $A \o H$, and the discussion there shows that the Radford biproduct $A \o H$ is cosemisimple if and only if its two factors~$A$ and~$H$ are cosemisimple. In the second step of this proof, we have also described the integrals in the Radford biproduct, and this description together with the first assertion of the proposition and its converse shows that the Radford biproduct $A \o H$ is semisimple if and only if its two factors~$A$ and~$H$ are semisimple. Analogous statements hold for~$B$ and~$H'$.

Now it also follows from the discussion in the second and third step of the preceding proof that in characteristic zero a Hopf algebra is semisimple if and only if it is finite-dimensional and cosemisimple. Since both~$H$ and~$H'$ are semisimple by assumption, and therefore finite-dimensional and cosemisimple, we see first that~$A$, which is also semisimple by assumption, is finite-dimensional and cosemisimple, and second that~$B$ is semisimple if and only if it is cosemisimple.

\item
Because~$B$ is a subcoalgebra of~$A$, the restriction of the integral \mbox{$\lma \in A^*$} described above to~$B$ is still an integral on~$B$. Because $\A \in B$ by assumption and
$\lma(\A)=\dim(A) \neq 0$, this restriction is nonzero, and~$B$ is cosemisimple by the version of Maschke's theorem for Yetter-Drinfel'd Hopf algebras in the third step of the preceding proof that we have already mentioned several times. As we saw in the first step, this implies that~$B$ is semisimple. 

\item
We now consider the corresponding integrals $\lmb \in B^*$ and~$\Lmb \in B$ for~$B$, subject to the analogous normalization conditions $\eb(\Lmb)=1$ and~$\lmb(\Lmb)=1$. Because the spaces of integrals are one-dimensional (cf.~\cite{Doi2}, Thm.~3, p.~3066; \cite{FMS}, Cor.~5.8, p.~4885; \cite{SoRib}, Prop.~2.10, p.~432), the restriction of the integral~$\lma \in A^*$ to~$B$ must be proportional to~$\lmb$. Because $\lma(\A)=\dim(A)$ and \mbox{$\lmb(\A)=\dim(B)$} by the preceding proposition, we must have
$$\lma |_B = \frac{\dim(A)}{\dim(B)} \lmb$$
which clearly implies
$$\lma(\Lmb) = \frac{\dim(A)}{\dim(B)}$$
Since $\Lmb$ is an idempotent and $\lma$ is, according to the preceding proposition, the character of the (left or right) regular representation, the left-hand side of this equation is an integer. This shows that $\dim(B)$ divides~$\dim(A)$. 

\item
If $l:=\dim(A)/\dim(B)$, the above equality states that the restriction of the character of the (left or right) regular representation of~$A$ to~$B$ is~$l$ times the character of the (left or right) regular representation of~$B$. This implies that
$$A \cong B^l$$
as (left or right) $B$-modules; in other words, $A$ is free of rank~$l$ over $B$, as asserted. 
\qed
\end{pflist}
\end{pf}

It should be mentioned that there are other freeness results for Yetter-Drinfel'd Hopf algebras in the literature in which the smaller Yetter-Drinfel'd Hopf algebra is not necessarily defined over the same Hopf algebra as the larger Yetter-Drinfel'd Hopf algebra: Such results for Nichols algebras, which are certain special Yetter-Drinfel'd Hopf algebras, have been obtained by M.~Gra\~{n}a (cf.~\cite{Grana}, Thm.~3.8, p.~250).

\section{Finite Abelian Groups} \label{Sec:FinAbelGroups}
\subsection[Finite abelian groups]{} \label{FinAbelGroups}
We now specialize the situation to the case where $H=K[G]$ is the group ring of a finite abelian group over an algebraically closed field~$K$ of characteristic zero. A left Yetter-Drinfel'd module~$V$ over~$H$ is in particular a module over~$H$, and this module structure can be described equivalently by a representation
$$G \rightarrow \Aut(V),~g \mapsto \phi_g$$
of~$G$. $V$ is also a left comodule over~$H$, i.e., a $G$-graded vector space (cf.~\cite{M}, Ex.~1.6.7, p.~12), and as the Yetter-Drinfel'd condition reduces to the dimodule condition in this setting, each~$\phi_g$ preserves the homogeneous components. However, we will use an equivalent, but different description of this comodule structure. As we discussed already in Paragraph~\ref{DualHopf}, every left $H$-comodule structure leads to a right module structure over~$H^*$. Since the character group $\hat{G} := \Hom(G,K^\times)$ of group homomorphisms from~$G$ to~$K^\times$ becomes a basis of~$H^*$ if we extend group homomorphisms $\gamma: G \rightarrow K^\times$ linearly to all of~$H$, we can consider~$H^*$ as the group ring~$K[\hat{G}]$ of the character group. This means that we also have a representation 
$$\hat{G} \rightarrow \Aut(V),~\gamma \mapsto \psi_\gamma$$
of~$\hat{G}$. Because $\hat{G}$ is abelian, the fact that this representation is derived from a right $H^*$-module structure, and not from a left $H^*$-module structure, does not pose a problem. Explicitly, $\psi_\gamma$ is then given by the formula $\psi_\gamma(v) = \gamma(v^\1) v^\2$. The Yetter-Drinfel'd condition implies that these two actions commute, i.e., that we have $\psi_\gamma \circ \phi_g = \phi_g \circ \psi_\gamma$ for all $g \in G$ and all $\gamma \in \hat{G}$, because~$\phi_g$ preserves the homogeneous components.

Since we have just explained that an $H$-comodule structure is the same as a representation of~$\hat{G}$, it must be possible to compute the coaction from this representation. This is achieved by the following formula:
\begin{prop}
For all $v \in V$, we have
$$\delta(v) = \frac{1}{|G|} \sum_{\gamma \in \hat{G}} \sum_{g \in G} \gamma(g^{-1}) \; g \o \psi_\gamma(v)$$
If~$W$ is a second Yetter-Drinfel'd module and the corresponding representations for~$W$ are also denoted by~$\phi$ and~$\psi$, the quasisymmetry and its inverse have the form
\begin{align*}
&\sigma_{V,W}(v \o w) = 
\frac{1}{|G|} \sum_{\gamma \in \hat{G}} \sum_{g \in G} \gamma(g^{-1}) \; \phi_g(w) \o \psi_\gamma(v) \\
&\sigma^{-1}_{V,W}(w \o v) = 
\frac{1}{|G|} \sum_{\gamma \in \hat{G}} \sum_{g \in G} \gamma(g) \; \psi_\gamma(v) \o \phi_g(w)
\end{align*}
\end{prop}
\begin{pf}
The orthogonality relations for finite abelian groups assert that
$$\sum_{g \in G} \gamma(g^{-1}) \gamma'(g) = |G| \delta_{\gamma,\gamma'}$$
for all $\gamma, \gamma' \in \hat{G}$. Therefore, applying $\gamma' \o \id_V$ to the right-hand side of the first assertion about the form of~$\delta(v)$ yields
\begin{align*}
(\gamma' \o \id_V&) \left(\frac{1}{|G|} 
\sum_{\gamma \in \hat{G}} \sum_{g \in G} \gamma(g^{-1}) \; g \o \psi_\gamma(v)\right)
= \frac{1}{|G|} \sum_{\gamma \in \hat{G}} \sum_{g \in G} \gamma(g^{-1}) \gamma'(g) \psi_\gamma(v) \\
&= \sum_{\gamma \in \hat{G}} \delta_{\gamma,\gamma'} \psi_\gamma(v) = \psi_{\gamma'}(v)
= \gamma'(v^\1) v^\2 = (\gamma' \o \id_V)(\delta(v))
\end{align*}
Because the set of all $\gamma' \in \hat{G}$ is a basis of~$H^*$, this proves the first assertion about the form of~$\delta(v)$. The second assertion follows directly by inserting the first assertion into the form of the quasisymmetry and its inverse stated in Paragraph~\ref{LeftYetMod}. 
\qed
\end{pf}

As we discussed in Paragraph~\ref{RightYetMod}, right Yetter-Drinfel'd modules are just left Yetter-Drinfel'd modules over~$H^{\op \cop}$. Therefore, the corresponding formula for the coaction reads
$$\delta(v) = \frac{1}{|G|} \sum_{\gamma \in \hat{G}} \sum_{g \in G} \gamma(g^{-1}) \; \psi_\gamma(v) \o g$$
where $\psi_\gamma$ is analogously given as $\psi_\gamma(v) = v^\1 \gamma(v^\2)$.

The proposition above has the following consequence:
\begin{corollary}
For $v \in V$ and $Q \subset \hat{G}$, the following two conditions are equivalent: 
\begin{enumerate}
\item 
$\psi_\gamma(v) = v$ for all $\gamma \in Q$

\item
$\delta(v) \in K[Q^\perp] \o V$
\end{enumerate}
In this case, the coaction has the form
$$\delta(v) = \frac{1}{|G|} \sum_{\gamma \in \hat{G}} \sum_{g \in Q^\perp} \gamma(g^{-1}) \; g \o \psi_\gamma(v)$$
\end{corollary}
\begin{pf}
The first assertion of the preceding proposition can be written in the form
$$\delta(v) = \frac{1}{|G|}  \sum_{g \in G} g \o 
\left(\sum_{\gamma \in \hat{G}} \gamma(g^{-1}) \psi_\gamma(v)\right)$$
which shows that the homogeneous component of degree~$g$ of~$v$ is given by the formula
$\frac{1}{|G|} \sum_{\gamma \in \hat{G}} \gamma(g^{-1}) \psi_\gamma(v)$. To see that the first condition implies the second, we have to show that this component vanishes if $g \notin Q^\perp$. But then there exists $\gamma' \in Q$ with $\gamma'(g) \neq 1$. As we have 
\begin{align*}
\sum_{\gamma \in \hat{G}} \gamma(g^{-1}) \psi_\gamma(v) &=
\sum_{\gamma \in \hat{G}} \gamma(g^{-1}) \psi_\gamma(\psi_{\gamma'}(v)) \\
&= \sum_{\gamma \in \hat{G}} (\gamma \gamma'^{-1})(g^{-1}) \psi_\gamma(v) =
\gamma'(g) \sum_{\gamma \in \hat{G}} \gamma(g^{-1}) \psi_\gamma(v)
\end{align*}
this yields that $\sum_{\gamma \in \hat{G}} \gamma(g^{-1}) \psi_\gamma(v) = 0$, as asserted.
To see that the second condition implies the first, we observe that for $\gamma \in Q$ and $g \in Q^\perp$, we have $\gamma(g) = 1 = \eh(g)$, so that $\delta(v) \in K[Q^\perp] \o V$ entails
$$\psi_\gamma(v) = \gamma(v^\1) v^\2 = \eh(v^\1) v^\2 = v$$
as asserted. 
\qed
\end{pf}

Although we have not assumed in this corollary that~$Q$ is a subgroup of~$\hat{G}$, this point is not decisive, because~$v$ is stabilized by~$Q$ if and only if it is stabilized by the subgroup generated by~$Q$.

As in the case of the proposition above, there is a version of this corollary for a right Yetter-Drinfel'd module~$V$, which states in particular that 
$$\delta(v) = \frac{1}{|G|} \sum_{\gamma \in \hat{G}} \sum_{g \in Q^\perp} \gamma(g^{-1}) \; \psi_\gamma(v) \o g$$
if $v \in V$ satisfies $\psi_\gamma(v) = v$ for all $\gamma \in Q$.

For two left Yetter-Drinfel'd modules~$V$ and~$W$ and a subgroup~$Q \subset \hat{G}$ that stabilizes an element~$v \in V$, we can clearly insert the form of the coaction obtained in the preceding corollary into the formula for the quasisymmetry. However, we can do more and also take an analogous subgroup of~$G$ into account:
\begin{lemma}
If we assume that elements~$v \in V$ and~$w \in W$ are stabilized by subgroups~$T \subset G$ and~$Q \subset \hat{G}$ in the sense that $\psi_\gamma(v) = v$ for all $\gamma \in Q$ and
$\phi_g(w) = w$ for all $g \in T$, we have
\begin{align*}
&\sigma_{V,W}(v \o w) = 
\frac{1}{|Q^\perp||Q \cap T^\perp|} \sum_{\gamma \in T^\perp}  \sum_{g \in Q^\perp} 
\gamma(g^{-1}) \; \phi_g(w) \o \psi_{\gamma}(v) \\
&\sigma^{-1}_{V,W}(w \o v) = 
\frac{1}{|Q^\perp||Q \cap T^\perp|} \sum_{\gamma \in T^\perp} \sum_{g \in Q^\perp} \gamma(g) \; \psi_\gamma(v) \o \phi_g(w)
\end{align*}
\end{lemma}
\begin{pf}
By the preceding corollary, we have 
$$\delta(v) = 
\frac{1}{|G|} \sum_{\gamma \in \hat{G}} \sum_{g \in Q^\perp} 
\gamma(g^{-1}) \; g \o \psi_\gamma(v)$$
As in the proof of the preceding proposition, we insert this formula into the equations for the quasisymmetry and its inverse to get
\begin{align*}
&\sigma_{V,W}(v \o w) = 
\frac{1}{|G|} \sum_{\gamma \in \hat{G}} \sum_{g \in Q^\perp} 
\gamma(g^{-1}) \; \phi_g(w) \o \psi_\gamma(v) \\
&\sigma^{-1}_{V,W}(w \o v) = 
\frac{1}{|G|} \sum_{\gamma \in \hat{G}} \sum_{g \in Q^\perp} 
\gamma(g) \; \psi_\gamma(v) \o \phi_g(w)
\end{align*}
Now suppose that $\gamma \notin Q T^\perp = (Q^\perp \cap T)^\perp$. Then there exists 
$h \in Q^\perp \cap T$ such that $\gamma(h) \neq 1$, and since
\begin{align*}
\sum_{g \in Q^\perp} \gamma(g)  \phi_g(w) &= 
\sum_{g \in Q^\perp} \gamma(gh)  \phi_g(\phi_h(w)) 
= \gamma(h) \sum_{g \in Q^\perp} \gamma(g)  \phi_g(w)
\end{align*}
we must have $\sum_{g \in Q^\perp} \gamma(g) \phi_g(w) = 0$. Applying this argument to~$\gamma^{-1}$ instead of~$\gamma$, we also get that 
$\sum_{g \in Q^\perp} \gamma(g^{-1}) \phi_g(w) = 0$.
The formulas above therefore reduce to 
\begin{align*}
&\sigma_{V,W}(v \o w)  = 
\frac{1}{|G|} \sum_{\gamma \in Q T^\perp} \sum_{g \in Q^\perp} 
\gamma(g^{-1}) \; \phi_g(w) \o \psi_\gamma(v)\\
&\sigma^{-1}_{V,W}(w \o v) = 
\frac{1}{|G|} \sum_{\gamma \in Q T^\perp} \sum_{g \in Q^\perp} 
\gamma(g) \; \psi_\gamma(v) \o \phi_g(w)
\end{align*}
Now the map
$T^\perp \times Q \rightarrow Q T^\perp,~(\gamma, \gamma') \mapsto \gamma \gamma'$
is a surjective group homomorphism whose kernel consists of the pairs 
$(\gamma, \gamma^{-1})$ with $\gamma \in Q \cap T^\perp$. We can therefore write the formula for the quasisymmetry in the form
\begin{align*}
\sigma_{V,W}(v \o w) &= 
\frac{1}{|G||Q \cap T^\perp|} 
\sum_{\gamma \in  T^\perp} \sum_{\gamma' \in Q} \sum_{g \in Q^\perp} 
(\gamma \gamma')(g^{-1}) \; \phi_g(w) \o \psi_{\gamma\gamma'}(v) \\
&= \frac{|Q|}{|G||Q \cap T^\perp|} 
\sum_{\gamma \in  T^\perp}  \sum_{g \in Q^\perp} 
\gamma(g^{-1}) \; \phi_g(w) \o \psi_{\gamma}(v)
\end{align*}
Since the formula for~$\sigma^{-1}_{V,W}(w \o v)$ can be rewritten in the same way, this proves the assertion, as we have $|G|=|Q||Q^\perp|$.
\qed
\end{pf}

The coefficient $|Q^\perp||Q \cap T^\perp|$ that appears in the equations stated in this lemma seems to be asymmetric in~$T$ and~$Q$. However, as we will discuss after Lemma~\ref{ChangeGroup}, this is indeed not the case.

In the case of two right Yetter-Drinfel'd modules~$V$ and~$W$ and elements~$v \in V$ and~$w \in W$ that are stabilized by subgroups~$T \subset G$ and~$Q \subset \hat{G}$ in the sense that $\phi_g(v) = v$ for all $g \in T$ and $\psi_\gamma(w) = w$ for all $\gamma \in Q$, where~$\phi_g$ now denotes the right action of~$g$, the corresponding formulas read
\begin{align*}
&\sigma_{V,W}(v \o w) = 
\frac{1}{|Q^\perp||Q \cap T^\perp|} \sum_{\gamma \in T^\perp} \sum_{g \in Q^\perp} 
\gamma(g^{-1}) \; \psi_{\gamma}(w) \o \phi_g(v) \\
&\sigma^{-1}_{V,W}(w \o v) = 
\frac{1}{|Q^\perp||Q \cap T^\perp|} \sum_{\gamma \in T^\perp} \sum_{g \in Q^\perp} 
\gamma(g) \; \phi_g(v) \o \psi_\gamma(w)
\end{align*}
These equations follow from the preceding lemma by considering right Yetter-Drinfel'd modules as left Yetter-Drinfel'd modules over~$H^{\op\cop}$, keeping in mind that we then have to reverse the order of the tensor factors, as explained in Paragraph~\ref{RightYetMod}. Alternatively, one can of course also adapt the above proof directly to this situation, using the formula for the right coaction of~$H$ given after the preceding corollary.

It is important to note that the situation is entirely symmetric in~$\phi_g$ and~$\psi_\gamma$: As we have seen in Paragraph~\ref{DualHopf}, we can view a left Yetter-Drinfel'd module~$V$ as a right Yetter-Drinfel'd module over~$H^*\cong K[\hat{G}]$. Our discussion at the beginning of this paragraph already shows that the representation of~$\hat{G}$ that is analogous to the representation $g \mapsto \phi_g$ of~$G$ in this setup is exactly~$\gamma \mapsto \psi_\gamma$. But the compatibility conditions in Paragraph~\ref{DualHopf} also show that the $H^*$-coaction, when turned into an action of~$H^{**} \cong H$ and then into a representation of~$G$, leads exactly to the representation $g \mapsto \phi_g$. Under this operation, the roles of the groups~$T$ and~$Q$ are reversed. In view of the symmetry property of the coefficient that we have already mentioned and will establish after Lemma~\ref{ChangeGroup}, the formulas for the quasisymmetry and its inverse in the case of right Yetter-Drinfel'd modules therefore also arise from the corresponding formulas in the case of left Yetter-Drinfel'd modules given in the preceding lemma by interchanging~$\phi_g$ and~$\psi_\gamma$.

\subsection[Change of groups]{} \label{ChangeGroup}
We now turn our attention from Yetter-Drinfel'd modules to Yetter-Drin\-fel'd Hopf algebras, still over the group ring of our finite abelian group~$G$. So, we consider now a left Yetter-Drinfel'd Hopf algebra~$A$ over \mbox{$H=K[G]$} and suppose that~$T \subset G$ and $Q \subset \hat{G}$ are subgroups that act trivially on~$A$, so that we have $\phi_g = \psi_\gamma = \id_A$ for all $g \in T$ and all $\gamma \in Q$. We define the group
$$G':= Q^\perp/(T \cap Q^\perp)$$
Its character group can be described as follows:
\begin{lemma}
$\hat{G}' \cong T^\perp/(Q \cap T^\perp)$
\end{lemma}
\begin{pf}
Every character $\gamma \in T^\perp$ induces a character
$$G' \rightarrow K^\times,~\bar{g} \mapsto \gamma(g)$$
so that we get a group homomorphism from~$T^\perp$ to~$\hat{G}'$ whose kernel is $Q \cap T^\perp$. Factoring over the kernel, we get an injective homomorphism from \mbox{$T^\perp/(Q \cap T^\perp)$} to~$\hat{G}'$. This shows that $|T^\perp/(Q \cap T^\perp)| \le |\hat{G}'|$, and it also shows that this map will be an isomorphism if $|\hat{G}'| = |T^\perp/(Q \cap T^\perp)|$.

Now every $g \in Q^\perp$ also induces a character
$$T^\perp/(Q \cap T^\perp) \rightarrow K^\times,~\bar{\gamma} \mapsto \gamma(g)$$
so that we get a group homomorphism from~$Q^\perp$ to the character group of 
$T^\perp/(Q \cap T^\perp)$ whose kernel is $T \cap Q^\perp$. Factoring over the kernel, we get 
an injective homomorphism from $G'= Q^\perp/(T \cap Q^\perp)$ to the character group of 
$T^\perp/(Q \cap T^\perp)$, which shows that $|\hat{G}'| = |G'| \le |T^\perp/(Q \cap T^\perp)|$.
\qed
\end{pf}

We note that this lemma implies that 
$$|Q^\perp|/|T \cap Q^\perp| = |G'|=|\hat{G}'| = |T^\perp|/|Q \cap T^\perp|$$
so that $|Q^\perp||Q \cap T^\perp| = |T^\perp||T \cap Q^\perp|$, which means that the coefficient appearing in Lemma~\ref{FinAbelGroups} is indeed symmetric in~$T$ and~$Q$, as we claimed there.

Clearly, the action of~$G$ on~$A$ can be restricted to~$Q^\perp$, and this action factors over~$G'$, so that~$A$ becomes a left $K[G']$-module. Similarly, the action of~$\hat{G}$ on~$A$ can be restricted to~$T^\perp$, and this action factors over~$T^\perp/(Q \cap T^\perp)$, which in view of the preceding lemma yields an action of~$\hat{G}'$ on~$A$. This action can then again be used to introduce a left $K[G']$\nobreakdash-comodule structure on~$A$, which by Proposition~\ref{FinAbelGroups} has the explicit form
\begin{align*}
\delta(a) &= 
\frac{1}{|G'|} \sum_{\gamma \in \hat{G}'} \sum_{g \in G'} 
\gamma(g^{-1}) \; g \o \psi_\gamma(a) 
\end{align*}
As it turns out, these structures fit together:
\begin{prop}
$A$ is a Yetter-Drinfel'd Hopf algebra over~$K[G']$.
\end{prop}
\begin{pf}
The fact that the two actions commute means that the Yetter-Drinfel'd compatibility condition is satisfied. Furthermore, both $G'$ and~$\hat{G}'$ act by algebra and coalgebra homomorphisms, so that~$A$ is a Yetter-Drinfel'd algebra and a Yetter-Drinfel'd coalgebra over~$K[G']$. The antipode~$\sa$ is also linear and colinear over~$K[G']$ and still satisfies the antipode axioms. To establish the assertion, it therefore remains to show that the quasisymmetry that comes from the $K[G']$-structure, which we denote by~$\sigma'_{A,A}$, coincides with the quasisymmetry~$\sigma_{A,A}$ that comes from the $K[G]$-structure. As we just saw, the coaction is given by
\begin{align*}
\delta(a) &= 
\frac{1}{|G'|} \sum_{\gamma \in \hat{G}'} \sum_{g \in G'} 
\gamma(g^{-1}) \; g \o \psi_\gamma(a) \\
&=  \frac{|T \cap Q^\perp|}{|Q^\perp|} \sum_{\bar{\gamma} \in T^\perp/(Q \cap T^\perp)} \sum_{\bar{g} \in Q^\perp/(T \cap Q^\perp)} \gamma(g^{-1}) \;
\bar{g} \o \psi_{\bar{\gamma}}(a) \\
&=  \frac{1}{|Q^\perp|} \sum_{\bar{\gamma} \in T^\perp/(Q \cap T^\perp)} \sum_{g \in Q^\perp} \gamma(g^{-1}) \; \bar{g} \o \psi_{\bar{\gamma}}(a)   \\
&=  \frac{1}{|Q^\perp||Q \cap T^\perp|} \sum_{\gamma \in T^\perp} \sum_{g \in Q^\perp} \gamma(g^{-1}) \; \bar{g} \o \psi_{\gamma}(a)  
\end{align*}
Using Heyneman-Sweedler notation with respect to the $K[G']$-structure, we find the expression
$$\sigma'_{A,A}(a \o b) = (a^\1.b) \o a^\2
=  \frac{1}{|Q^\perp||Q \cap T^\perp|} \sum_{\gamma \in T^\perp} \sum_{g \in Q^\perp} \gamma(g^{-1}) \phi_g(b) \o \psi_{\gamma}(a)$$
for the quasisymmetry. It therefore follows from Lemma~\ref{FinAbelGroups} that 
$$\sigma'_{A,A}(a \o b) = \sigma_{A,A}(a \o b)$$
for all $a,b \in A$. This implies that the algebra structure determined by~$\sigma'_{A,A}$ on the second tensor power of~$A$ coincides with the algebra structure of $A \hat{\o} A$, which is what we needed to show. 
\qed
\end{pf}

\section{The Commutative Case} \label{Sec:CommCase}
\subsection[Basic notions]{} \label{BasicNot}
In this section, we turn to the main object of our study: We consider a commutative semisimple left Yetter-Drinfel'd Hopf algebra~$A$ over the group ring~$H=K[G]$ of a finite abelian group~$G$. As in Section~\ref{Sec:FinAbelGroups}, the base field~$K$ is assumed to be algebraically closed of characteristic zero, except for a minor comment in Paragraph~\ref{TrivThm}. Then~$A$ is finite-dimensional by Proposition~\ref{Int}, and because~$K$ is algebraically closed, semisimplicity implies that $A$~has a basis consisting of primitive idempotents, any two of which are orthogonal. The set of all primitive idempotents of~$A$ will be denoted by~$\E$. Since the mappings~$\phi_g$ and~$\psi_\gamma$ defined in Paragraph~\ref{FinAbelGroups} are algebra automorphisms, they permute these primitive idempotents. The dual basis of~$\E$ consists exactly of the \mbox{(one-dimensional)} characters, i.e., the algebra homomorphisms from~$A$ to~$K$, and the set of all these characters is precisely the set~$G(A^*)$ of all group-like elements in the dual~$A^*$. 

\begin{samepage}
With each primitive idempotent, we associate a number of objects that will play an important role in the sequel:
\begin{defn}
Suppose that~$e \in \E$ is a primitive idempotent of~$A$.
\begin{enumerate}
\item 
The character associated with~$e$ is the unique algebra homomorphism to the base field $\eta_e: A \rightarrow K$ with the property that $\eta_e(e) = 1$ and $\eta_e(e') = 0$ for all primitive idempotents $e' \neq e$. 

\item 
The inertia group of~$e$ is the group $T_e :=\{g \in G \mid \phi_g(e) = e\}$. 

\item 
The isotropy group of~$e$ is the group $Q_e :=\{\gamma \in \hat{G} \mid \psi_\gamma(e) = e\}$. 

\item 
The index group of~$e$ is the group $G_e := Q_e^\perp/(T_e \cap Q_e^\perp)$. 

\item 
The index of~$e$ is the number $|G_e|$. 

\item 
The orbit of~$e$ under the action of~$Q_e^\perp$ is the set $O_e := \{\phi_g(e) \mid g \in Q_e^\perp\}$.

\item 
The full orbit of~$e$ is the set 
$\hat{O}_e := \{\psi_\gamma(\phi_g(e)) \mid g \in G, \gamma \in \hat{G}\}$.

\item 
The stability set of~$e$ is the set 
$\check{O}_e := \{e' \in \E \mid T_e \subset T_{e'}\; \text{and} \; Q_e \subset Q_{e'}\}$.
\end{enumerate}
\end{defn}
\end{samepage}
The ideals spanned by these sets will be denoted by
$$I_e := \Span(O_e) \qquad \qquad \hat{I}_e := \Span(\hat{O}_e) \qquad \qquad 
\check{I}_e := \Span(\check{O}_e)$$
and the unique ideals that complement them will be denoted by~$J_e$, $\hat{J}_e$, and~$\check{J}_e$, respectively. Stated differently, $J_e$, $\hat{J}_e$, and~$\check{J}_e$ are spanned by all primitive idempotents that are not contained in~$O_e$, $\hat{O}_e$, and $\check{O}_e$, respectively, and they satisfy
$A= I_e \oplus J_e = \hat{I}_e \oplus\hat{J}_e = \check{I}_e \oplus \check{J}_e$. The characters that correspond to these sets of idempotents will be denoted by
$$U_e := \{\eta_{e'} \mid e' \in O_e\} \quad \qquad 
\hat{U}_e := \{\eta_{e'} \mid e' \in \hat{O}_e\} \quad \qquad 
\check{U}_e := \{\eta_{e'} \mid e' \in \check{O}_e\}$$
so that
$$J_e^\perp = \Span(U_e) \qquad \qquad \hat{J}_e^\perp = \Span(\hat{U}_e) \qquad \qquad
\check{J}_e^\perp = \Span(\check{U}_e)$$
By construction, we have $O_e \subset \hat{O}_e \subset \check{O}_e$ and therefore
$$I_e \subset \hat{I}_e \subset \check{I}_e  \qquad \qquad
J_e \supset \hat{J}_e \supset \check{J}_e \qquad \qquad
U_e \subset \hat{U}_e \subset \check{U}_e$$
It follows directly from these definitions that $\hat{O}_e$ and~$\check{O}_e$, and accordingly also $\hat{I}_e$, $\check{I}_e$, $\hat{J}_e$, $\check{J}_e$, $\hat{U}_e$, and~$\check{U}_e$, are stable under both~$G$ and~$\hat{G}$. Also by construction we have $\dim(I_e)= |O_e| = |G_e|$.
We will furthermore speak of the index of the character~$\eta_e$, which is by definition the index of the corresponding primitive idempotent~$e$. These definitions should be compared with \cite{SoYp}, Def.~2.3, p.~27. 

\begin{samepage}
We can compute a primitive idempotent from its character via the integrals introduced in Paragraph~\ref{Int}:
\begin{prop}
If $e \in \E$ is a primitive idempotent, we have
\begin{align*}
e &= \eta_e^{-1}(\Lma_\1) \Lma_\2 = \eta_e(\Lma_\1) \sa(\Lma_\2) \\
&= \Lma_\1 \eta_e^{-1}(\Lma_\2) = \sa(\Lma_\1) \eta_e(\Lma_\2)
\end{align*}
Furthermore, we have $\lma(e)=1$.
\end{prop}
\end{samepage}
\begin{pf}
We know from Proposition~\ref{Int} that
$$e = \lma(e \Lma_\1) \sa(\Lma_\2) = \lma(e) \eta_e(\Lma_\1) \sa(\Lma_\2)$$
and also $e = \sa(\Lma_\1) \eta_e(\Lma_\2) \lma(e)$.
If we apply~$\eta_e$ to the first equation, we get that
$$1 = \eta_e(e) = \lma(e) \eta_e(\Lma_\1) \eta_e(\sa(\Lma_\2))= \lma(e) \ea(\Lma) \eta_e(\A) = \lma(e)$$
Inserting this back into these equations, we get two of the four forms of~$e$. The remaining two follow by substituting the other form of the Casimir element into the ones already established, as we have
$\eta_e^{-1} = \eta_e \circ \sa$.
\qed
\end{pf}

The fact that $\lma(e)=1$ gives another proof for the assertion that~$\lma$ is the character of the regular representation made in Proposition~\ref{Int}, since the character of the regular representation clearly also takes the value~$1$ on primitive idempotents, and the primitive idempotents form a basis, as discussed at the beginning of this paragraph. 

Let us note a remarkable consequence of this proposition:
\begin{corollary}
For primitive idempotents~$e$ and~$e'$, we have $\eta_e^{-1}(e') = \eta_{e'}^{-1}(e)$.
\end{corollary}
\begin{pf}
We have seen in the last proposition that 
$$e = \eta_e^{-1}(\Lma_\1) \Lma_\2 = \Lma_\1 \eta_e^{-1}(\Lma_\2)$$
which obviously implies that
$$\eta_{e'}^{-1}(e) = \eta_e^{-1}(\Lma_\1) \eta_{e'}^{-1}(\Lma_\2) 
= \eta_{e'}^{-1}(\Lma_\1) \eta_e^{-1}(\Lma_\2)$$
If we interchange~$e$ and~$e'$ in this equation, it becomes
$$\eta_e^{-1}(e') = \eta_{e'}^{-1}(\Lma_\1) \eta_e^{-1}(\Lma_\2) = \eta_e^{-1}(\Lma_\1) \eta_{e'}^{-1}(\Lma_\2)$$
The assertion then follows by comparing these two equations.~\qed
\end{pf}

We note that the proof of this corollary shows again that $\Lma$ is cocommutative, as we already stated in Proposition~\ref{Int}: We have just seen in this proof that
$$(\eta_e^{-1} \o \eta_{e'}^{-1})(\Lma_\1 \o \Lma_\2) = (\eta_e^{-1} \o \eta_{e'}^{-1})(\Lma_\2 \o \Lma_\1)$$
Since the tensors $\eta_e^{-1} \o \eta_{e'}^{-1}$ form a basis of~$A^* \o A^*$, this implies that we have
$\Lma_\1 \o \Lma_\2 = \Lma_\2 \o \Lma_\1$.

\subsection[Ideals in the tensor product]{} \label{IdealTens}
We now investigate the ideal structure of the algebra~$A \hat{\o} A$ introduced in Paragraph~\ref{YetDrinfHopf} in our situation. For this, we need to associate with two primitive idempotents~$e$ and~$e'$ two ideals which are defined in a similar way to the ideals~$I_e$ and~$I_{e'}$, but take the interaction of~$e$ and~$e'$ into account. So, we consider the orbits
$$O := \{\phi_g(e) \mid g \in Q_{e'}^\perp\} \qquad \qquad
O' := \{\psi_\gamma(e') \mid \gamma \in T_{e}^\perp\}$$
and denote the ideals that they span by $I:=\Span(O)$ and $I':=\Span(O')$, respectively. Inside these ideals, we need to consider certain elements which arise from the idempotents by a discrete Fourier transform:
\begin{defb}
\begin{enumerate}
\item 
For every $\gamma \in \hat{G}$, we define 
$$w_\gamma := \frac{1}{|T_e \cap Q_{e'}^\perp|} \sum_{g \in Q_{e'}^\perp} \gamma(g^{-1}) \phi_g(e) \in I$$

\item
For every $g \in G$, we define
$$u_g := \frac{1}{|Q_{e'} \cap T_e^\perp|} 
\sum_{\gamma \in T_e^\perp} \gamma(g^{-1}) \psi_\gamma(e') \in I'$$
\end{enumerate}
\end{defb}

We record that in the case where $\gamma = \eh$ is the counit, $w_\gamma$ is just the sum of all elements of~$O$, which we denote by~$e_I$ (cf. \cite{SoYp}, Par.~6.6, p.~96). If we consider the ideal~$I$ as a non-unital subalgebra of~$A$, then~$e_I$ is the unit element of this algebra.

The following lemma lists the basic properties of the elements~$w_\gamma$:
\begin{lemb}
\begin{enumerate}
\item 
$O$ and $O'$ have the same cardinality~$m$.

\item 
For all $g \in Q_{e'}^\perp$, we have $\phi_g(w_\gamma) = \gamma(g) w_\gamma$.

\item
If $\gamma$ and~$\gamma'$ belong to the same coset in $\hat{G}/Q_{e'}$, we have $w_\gamma = w_{\gamma'}$.

\item
If $\gamma \notin T_e^\perp Q_{e'}$, we have $w_\gamma =  0$.

\item
If $\gamma \in T_e^\perp Q_{e'}$, we have $\eta_e(w_\gamma) = 1$.

\item
If $g_1,\ldots,g_m \in Q_{e'}^\perp$ is a system of representatives for the cosets in the quotient 
$Q_{e'}^\perp/(T_e \cap Q_{e'}^\perp)$, then we have 
$$w_\gamma = \sum_{i=1}^m \gamma(g_i^{-1}) \phi_{g_i}(e)$$
for all $\gamma \in T_e^\perp Q_{e'}$.

\item
If $\gamma_1,\ldots,\gamma_m \in T_e^\perp$ is a system of representatives
for the cosets in the quotient $T_e^\perp/(Q_{e'} \cap T_e^\perp)$, then
the elements $w_{\gamma_1},\ldots,w_{\gamma_m}$ form a basis of~$I$.

\item
For $g \in Q_{e'}^\perp$, we then have
$$\phi_g(e) = \frac{1}{m} \sum_{j=1}^m \gamma_j(g) w_{\gamma_j}$$
\end{enumerate}
\end{lemb}
\begin{pf}
\begin{pflist}
\item
We have $|O|=|Q_{e'}^\perp/(T_e \cap Q_{e'}^\perp)|$ and 
$|O'|=|T_e^\perp/(Q_{e'} \cap T_e^\perp)|$. Since we have seen in the proof of Lemma~\ref{ChangeGroup} that the cardinalities of these two factor groups are equal, the first assertion holds.

\item
The second and the third assertion follow more or less directly from the definition. For the fourth assertion, we note that if $\gamma \notin T_e^\perp Q_{e'} = (T_e \cap Q_{e'}^\perp)^\perp$, we can find $g \in T_e \cap Q_{e'}^\perp$ satisfying 
$\gamma(g) \neq 1$, which implies that
$$w_\gamma = \phi_g(w_\gamma) = \gamma(g) w_\gamma$$ 
by the second property, so that $w_\gamma = 0$.

\item
The fifth assertion holds since
\begin{align*}
\eta_e(w_\gamma) &= \frac{1}{|T_e \cap Q_{e'}^\perp|} \sum_{g \in Q_{e'}^\perp} \gamma(g^{-1}) \eta_e(\phi_g(e)) \\
&= \frac{1}{|T_e \cap Q_{e'}^\perp|} \sum_{g \in T_e \cap Q_{e'}^\perp} \gamma(g^{-1}) \eta_e(e) = 1
\end{align*}
where we have used in the last step that $\gamma(g)=1$ if $g \in T_e \cap Q_{e'}^\perp$ and 
$\gamma \in T_e^\perp Q_{e'}$.

\item
To prove the sixth property, we use the system of representatives to write the definition of~$w_\gamma$ in the form
$$w_\gamma = \frac{1}{|T_e \cap Q_{e'}^\perp|} \sum_{i=1}^m 
\sum_{g \in T_e \cap Q_{e'}^\perp} \gamma(g^{-1} g_i^{-1}) \phi_{g_i g}(e)$$
As in the previous step, we have for $g \in T_e \cap Q_{e'}^\perp$ and~$\gamma \in T_e^\perp Q_{e'}$ that $\phi_{g}(e) = e$ and $\gamma(g^{-1}) = 1$. Therefore, the summation over~$g$ can be carried out, which gives the assertion.

\item
For the seventh property, we have that $O=\{\phi_{g_1}(e),\ldots, \phi_{g_m}(e)\}$ is a basis of~$I$ by definition, and the base change 
$$w_{\gamma_j} = \sum_{i=1}^m \gamma_j(g_i^{-1}) \phi_{g_i}(e)$$
is given by an invertible matrix, as we saw in the proof of Lemma~\ref{ChangeGroup}. By the orthogonality relations for finite abelian groups, which we have already recalled in Paragraph~\ref{FinAbelGroups}, the inverse of this matrix is $(\gamma_j(g_i)/m)$, so that we get
$$\phi_{g_i}(e) = \frac{1}{m} \sum_{j=1}^m \gamma_j(g_i) w_{\gamma_j}$$
If $i \le m$ is chosen so that~$g$ and~$g_i$ are in the same coset of the quotient 
$Q_{e'}^\perp/(T_e \cap Q_{e'}^\perp)$, we can replace~$g_i$ by~$g$ on the left-hand and the right-hand side to get the last assertion.
\qed
\end{pflist}
\end{pf}

It follows from the discussion in Paragraph~\ref{DualHopf} together with Lemma~\ref{YetDrinfHopf} that we can consider \mbox{$A^{\op \cop} = A^{\cop}$} as a left Yetter-Drinfel'd Hopf algebra over the dual Hopf algebra $H^* = H^{* \op \cop} \cong K[\hat{G}]$. As we have explained at the end of Paragraph~\ref{FinAbelGroups}, the roles of~$\phi_g$ and~$\psi_\gamma$ are interchanged from this viewpoint. If we also interchange the idempotents~$e$ and~$e'$, then the elements~$w_\gamma$ turn into the elements~$u_g$. In this way, the following properties of the elements~$u_g$ follow directly from the corresponding properties of the elements~$w_\gamma$:
\begin{corb}
\begin{enumerate}
\item 
For all $\gamma \in T_e^\perp$, we have $\psi_\gamma(u_g) = \gamma(g) u_g$.

\item 
If $g$ and $g'$ belong to the same coset in~$G/T_e$, we have $u_g = u_{g'}$.

\item 
If $g \notin Q_{e'}^\perp T_e$, we have $u_g=0$.

\item
If $g \in Q_{e'}^\perp T_e$, we have $\eta_{e'}(u_g) = 1$.

\item
If $\gamma_1,\ldots,\gamma_m \in T_e^\perp$ is a system of representatives
for the cosets in the quotient $T_e^\perp/(Q_{e'} \cap T_e^\perp)$, then we have
$$u_g = \sum_{j=1}^m \gamma_j(g^{-1}) \psi_{\gamma_j}(e')$$
for all 
$g \in Q_{e'}^\perp T_e$.

\item
If $g_1,\ldots,g_m \in Q_{e'}^\perp$ is a system of representatives for the cosets in the quotient 
$Q_{e'}^\perp/(T_e \cap Q_{e'}^\perp)$, then the elements $u_{g_1}, \ldots, u_{g_m}$ form a basis of~$I'$.

\item
For $\gamma \in T_e^\perp$, we then have
$$\psi_\gamma(e') = \frac{1}{m} \sum_{i=1}^m \gamma(g_i) u_{g_i}$$
\end{enumerate} 
\end{corb}

\begin{samepage}
As in \cite{SoYp}, Prop.~6.5, p.~94, we can use these elements to describe the ideal structure of the algebra~$A \hat{\o} A$ introduced in Paragraph~\ref{YetDrinfHopf}:
\begin{propb}
\begin{enumerate}
\item
For all $a \in A$, we have $a^\1 \eta_{e'}(a^\2) \in K[Q_{e'}^\perp]$.

\item
For all $a \in I$ and all $g \in G$, we have $(u_g^\1.a) \o u_g^\2 = \phi_g(a) \o u_g$.

\item 
$I \hat{\o} K e'$ is a minimal left ideal of $A \hat{\o} A$.

\item
$A$ is an $A \hat{\o} A$-module with respect to the module structure
$$(a \o a').b = a (a'{}^\1.b) \eta_{e'}(a'{}^\2)$$

\item
For this module structure, $I$ is a simple $A \hat{\o} A$-submodule of~$A$.

\item
$I \hat{\o} I'$ is a minimal two-sided ideal of $A \hat{\o} A$.
\end{enumerate}
\end{propb}
\end{samepage}
\begin{pf}
\begin{pflist}
\item
The ideal~$\hat{I}_{e'}$ and its complement~$\hat{J}_{e'}$ are invariant under~$\hat{G}$ and therefore are \mbox{$H^*$-submodules}, i.e., $H$-subcomodules. If $a \in \hat{I}_{e'}$, we obtain from Corollary~\ref{FinAbelGroups} that $\delta(a) \in K[Q_{e'}^\perp] \o A$ and therefore 
$a^\1 \eta_{e'}(a^\2) \in K[Q_{e'}^\perp]$. On the other hand, if $a \in \hat{J}_{e'}$, then we have that 
$\delta(a) \in H \o \hat{J}_{e'}$ and therefore $a^\1 \eta_{e'}(a^\2) = 0 \in K[Q_{e'}^\perp]$.

\item
To prove the second assertion, we can assume that $a=w_\gamma$ for some $\gamma \in T_e^\perp$, because these elements span~$I$ by the preceding lemma. If $g \notin Q_{e'}^\perp T_e$, we have $u_g=0$ by the preceding corollary, so that the assertion is correct. If $g \in Q_{e'}^\perp T_e$, we also know from the preceding corollary that we can modify~$g$ so that $g \in Q_{e'}^\perp$, and the right-hand side also does not change under this modification. Since we know from Corollary~\ref{FinAbelGroups} that $\delta(u_g) \in K[Q_{e'}^\perp] \o A$, it follows from the preceding lemma that 
\begin{align*}
(u_g^\1.w_\gamma) \o u_g^\2 &= \gamma(u_g^\1) w_\gamma \o u_g^\2 = w_\gamma \o \psi_\gamma(u_g) \\
&= w_\gamma \o \gamma(g) u_g = \gamma(g)  w_\gamma \o u_g = \phi_g(w_\gamma) \o u_g
\end{align*}
establishing our claim.

\item
Now to see that $I \hat{\o} K e'$ is a left ideal, we observe that for $b \in I$ we have
$$(a \o a')(b \o e') = a (a'{}^\1.b) \o a'{}^\2 e' 
= a (a'{}^\1.b) \eta_{e'}(a'{}^\2) \o e'$$
and $(a'{}^\1.b) \eta_{e'}(a'{}^\2) \in I$ by the first assertion. The minimality will be discussed in the sixth step below.

\item
The fourth assertion can be proved directly by verifying the \mbox{module axioms}. However, we give a different proof with the help of an induced module construction: We can view~$A$ as a subalgebra of $A \hat{\o} A$ by mapping~$a$ to~$\A \o a$. Turning the base field~$K$ into an 
\mbox{$A$-module} by using~$\eta_{e'}$, we can then form the induced module $(A \hat{\o} A) \o_A K$, which is isomorphic to~$A$ as a vector space via
$b \mapsto (b \o \A) \o_A 1$. Now we have
\begin{align*}
(a \o a'&).((b \o \A) \o_A 1) = (a (a'{}^\1.b) \o a'{}^\2) \o_A 1 \\
&= (a (a'{}^\1.b) \o \A) \o_A a'{}^\2.1 = (a (a'{}^\1.b) \eta_{e'}(a'{}^\2) \o \A) \o_A 1
\end{align*}
If we use the isomorphism above to transfer this module structure to~$A$, we get the module structure that appears in the fourth assertion.

\item
The fact that~$I$ is an $A \hat{\o} A$-submodule of~$A$ with respect to this module structure follows from the first assertion. To see that~$I$ is simple, we consider a nonzero submodule~$U$ of~$I$. Because 
$(a \o \A).b = ab$, $U$ is in particular a left ideal of~$A$ and therefore spanned by a subset of the idempotents in~$O$. By the second assertion, we have for $g \in Q_{e'}^\perp$ and~$b \in I$ that  
$$(\A \o u_g).b = \phi_g(b) \eta_{e'}(u_g) = \phi_g(b)$$
where we have used the fact that $\eta_{e'}(u_g) = 1$ established in the preceding corollary. Therefore, the set of idempotents that span~$U$ is stable under~$Q_{e'}^\perp$. But since~$O$ is a $Q_{e'}^\perp$-orbit by definition, this proves $U=I$.

\item
Of course, the left ideal $I \hat{\o} K e'$ of $A \hat{\o} A$ is also a left $A \hat{\o} A$-module by multiplication. But our computation in the third step then shows that 
$$I \rightarrow I \hat{\o} K e',~a \mapsto a \o e'$$
is an isomorphism of $A \hat{\o} A$-modules. Therefore, the simplicity of~$I$ implies the minimality 
of~$I \hat{\o} K e'$.

\item
To see that $I \hat{\o} I'$ is a left ideal, we note that
$$I \hat{\o} I' = \bigoplus_{\tilde{e}' \in O'} I \hat{\o} K \tilde{e}'$$
Since the elements of~$O'$ have the same isotropy group~$Q_{e'}$ as~$e'$, the summands in this decomposition are left ideals by the third assertion. 

\item
To see that $I \hat{\o} I'$ is a right ideal, we have to show for $a \in I$, $a' \in I'$, and all $b,b' \in A$ that
$$(a \o a')(b \o b') = a (a'{}^\1.b) \o a'{}^\2 b' \in I \hat{\o} I'$$
For this, we introduce the ideal~$J$ spanned by all primitive idempotents that are not in~$O$. Then we have
$A = I \oplus J$, and both~$I$ and~$J$ are invariant under~$Q_{e'}^\perp$. From Corollary~\ref{FinAbelGroups}, we know that $\delta(a') \in K[Q_{e'}^\perp] \o A$. If $b \in J$, we therefore have $(a'{}^\1.b) \o a'{}^\2 \in J \o A$, so that
$$(a \o a')(b \o b') = a (a'{}^\1.b) \o a'{}^\2 b' = 0 \in I \hat{\o} I'$$
It remains to consider the case~$b \in I$. For this, we can assume that $a' = u_g$ for some $g \in Q_{e'}^\perp$, because these elements span~$I'$ by the preceding corollary. But then we get from the second assertion that
\begin{align*}
(a \o u_g)(b \o b') = a(u_g^\1.b) \o u_g^\2 b' = a \phi_g(b) \o u_g b' \in I \hat{\o} I'
\end{align*}
as asserted.

\item
The minimality of~$I \hat{\o} I'$ now follows from a result of W.~Burnside (cf.~\cite{FD}, Cor.~1.16, p.~45), which, since the argument in the fifth step also shows that~$I$ is simple as an $I \hat{\o} I'$-module, implies that the structure map $I \hat{\o} I' \rightarrow \End_K(I)$ is surjective, and therefore bijective by dimension considerations.
\qed
\end{pflist}
\end{pf}

It follows from this discussion that the left ideals $I \hat{\o} K \tilde{e}'$ for $\tilde{e}' \in O'$ that appear in the seventh step of the preceding proof are all isomorphic, because we have now seen that they are left ideals in a matrix algebra of the correct dimension. Equivalently, the $A \hat{\o} A$-module structure on~$I$ considered in the fifth assertion must be isomorphic to the one that arises in the same way from the character~$\eta_{\tilde{e}'}$ instead of~$\eta_{e'}$. If $\tilde{e}' = \psi_\gamma(e')$ for $\gamma \in T_e^\perp$, this isomorphism is indeed just given by multiplication with~$w_{\gamma^{-1}}$: Since $\eta_{\tilde{e}'} = \psi_{\gamma^{-1}}^*(\eta_{e'}) = \eta_{e'} \circ \psi_{\gamma^{-1}}$ is the character corresponding to~$\tilde{e}'$, we have
\begin{align*}
(a \o a'&).(b w_{\gamma^{-1}}) = a (a'{}^\1.(b w_{\gamma^{-1}})) \eta_{e'}(a'{}^\2)
= a (a'{}^\1.b) (a'{}^\2.w_{\gamma^{-1}}) \eta_{e'}(a'{}^\3) \\
&= a (a'{}^\1.b) (\gamma^{-1}(a'{}^\2) w_{\gamma^{-1}}) \eta_{e'}(a'{}^\3) 
= a (a'{}^\1.b) w_{\gamma^{-1}} \eta_{e'}(\psi_{\gamma^{-1}}(a'{}^\2)) \\
&= a (a'{}^\1.b) w_{\gamma^{-1}} \eta_{\tilde{e}'}(a'{}^\2)
\end{align*}
where the third equality uses the fact that $a^\1 \eta_{e'}(a^\2) \in K[Q_{e'}^\perp]$ established in the preceding proposition.

\subsection[Products of characters]{} \label{ProdChar}
For primitive idempotents~$e$ and~$e'$, we know from \cite{SoYp}, Prop.~1.5.2, p.~13 that the product~$\eta_e \eta_{e'}$ of the corresponding characters is again a character if $\eta_e$ is $H$-linear or $\eta_{e'}$ is colinear. Now we want to analyze what happens in general. We use the notation introduced in the preceding paragraph:
\begin{thmb}
\begin{enumerate}
\item
There are distinct characters $\omega_1,\ldots,\omega_m$ such that 
$$\eta_{\tilde{e}} \eta_{\tilde{e}'} \in \Span(\omega_1,\ldots,\omega_m)$$
for all $\tilde{e} \in O$ and all $\tilde{e}' \in O'$. 
In particular,
$\eta_e \eta_{e'} \in \Span(\omega_1,\ldots,\omega_m)$.

\item
If, for $\gamma \in T_e^\perp$, the orbit~$O$ is stable under $\psi_\gamma$, then
$\{\omega_1,\ldots,\omega_m \}$ is stable under $\psi_\gamma^*$.

\item
If, for $g \in Q_{e'}^\perp$, the orbit~$O'$ is stable under $\phi_g$, then
$\{\omega_1,\ldots,\omega_m \}$ is stable under $\phi_g^*$.
\end{enumerate} 
\end{thmb}
\begin{pf}
\begin{pflist}
\item
We have seen in Proposition~\ref{IdealTens} that~$I$ is a simple $A \hat{\o} A$-submodule of~$A$, with respect to the module structure described there. Since the coproduct $\da: A \rightarrow A \hat{\o} A$ is an algebra homomorphism, $I$ also becomes an $A$-module by restriction along~$\da$. Since $A$ is commutative, this $A$-module can be decomposed into a direct sum of one-dimensional $A$-modules. Therefore there is a basis $v_1,\ldots,v_m$ of~$I$ such that $a.v_k = \omega_k(a) v_k$ for characters $\omega_1,\ldots,\omega_m$. Explicitly, this equation asserts that
$$\omega_k(a) v_k = a.v_k = \da(a).v_k = a_\1 (a_\2{}^\1.v_k) \eta_{e'} (a_\2{}^\2)$$

\item 
For $\gamma \in T_e^\perp$, we can write~$w_\gamma$ in the form
$$w_\gamma = \sum_{k=1}^m \xi_k v_k$$ 
since $v_1,\ldots,v_m$ is a basis of~$I$. We then have 
\begin{align*}
a_\1 w_\gamma \psi^*_{\gamma}(\eta_{e'})(a_\2) 
&= a_\1 w_\gamma \gamma(a_\2{}^\1) \eta_{e'}(a_\2{}^\2) = a_\1 (a_\2{}^\1.w_\gamma) \eta_{e'}(a_\2{}^\2) \\
&= \sum_{k=1}^m \xi_k a_\1 (a_\2{}^\1.v_k) \eta_{e'}(a_\2{}^\2) 
=  \sum_{k=1}^m \xi_k \omega_k(a) v_k
\end{align*}
where we have used for the second equality again that $a^\1 \eta_{e'}(a^\2) \in K[Q_{e'}^\perp]$ by Proposition~\ref{IdealTens}. For $g \in Q_{e'}^\perp$, we have by Lemma~\ref{IdealTens} that
$$\phi_g^*(\eta_e)(w_\gamma) = \eta_e(\phi_g(w_\gamma)) = \gamma(g) \eta_e(w_\gamma) = \gamma(g)$$
Therefore we get by applying $\phi_g^*(\eta_e)$ to the equation above that
\begin{align*}
\gamma(g) \phi_g^*(\eta_e)(a_\1) \psi^*_{\gamma}(\eta_{e'})(a_\2) &= 
\phi_g^*(\eta_e)(a_\1 w_\gamma) \psi^*_{\gamma}(\eta_{e'})(a_\2) \\
&= \sum_{k=1}^m \xi_k \omega_k(a) \phi_g^*(\eta_e)(v_k)
\end{align*}
showing that
$$\phi_g^*(\eta_e) \psi^*_{\gamma}(\eta_{e'}) =
\gamma(g^{-1}) \sum_{k=1}^m \xi_k \phi_g^*(\eta_e)(v_k) \omega_k$$
If now $\tilde{e} \in O$ and $\tilde{e}' \in O'$, we have $\tilde{e} = \phi_{g^{-1}}(e)$ for some~$g \in Q_{e'}^\perp$ and similarly $\tilde{e}' = \psi_{\gamma^{-1}}(e')$ for some~$\gamma \in T_e^\perp$. The corresponding characters are $\eta_{\tilde{e}} = \phi_g^*(\eta_e)$ and 
$\eta_{\tilde{e}'} = \psi^*_{\gamma}(\eta_{e'})$. Therefore the last equation shows that 
$$\eta_{\tilde{e}} \eta_{\tilde{e}'} \in \Span(\omega_1,\ldots,\omega_m)$$
Since~$\eta_{e'}$ is invertible, the elements $\eta_{\tilde{e}} \eta_{e'}$ for $\tilde{e} \in O$ are linearly independent, which implies that the characters $\omega_1,\ldots,\omega_m$ must be distinct.

\item
To prove the second assertion, we argue as follows: We have
$$\Span(\omega_1,\ldots,\omega_m)
= \Span(\{\eta_{\tilde{e}} \eta_{\tilde{e}'} \mid \tilde{e} \in O, \tilde{e}' \in O'\})$$
By applying~$\psi_\gamma^*$, this becomes
$$\Span(\psi_\gamma^*(\omega_1),\ldots,\psi_\gamma^*(\omega_m))
= \Span(\{\psi_\gamma^*(\eta_{\tilde{e}}) \psi_\gamma^*(\eta_{\tilde{e}'}) 
\mid \tilde{e} \in O, \tilde{e}' \in O'\})$$
If $\gamma \in T_e^\perp$, then  $O'$ is stable under~$\psi_\gamma$ anyway. Therefore, if $O$ is stable under~$\psi_\gamma$, we get
$$\Span(\psi_\gamma^*(\omega_1),\ldots,\psi_\gamma^*(\omega_m)) =
\Span(\omega_1,\ldots,\omega_m)$$
which implies the assertion, since distinct group-like elements are linearly independent (cf.~\cite{M}, Def.~1.3.4, p.~4; \cite{RadfHopf}, Lem.~2.1.12, p.~25; \cite{Sw}, Prop.~3.2.1, p.~54). The proof of the third assertion is similar.
\qed
\end{pflist}
\end{pf}

The preceding theorem generalizes the results in \cite{SoYp}, Prop.~6.6, p.~95. It is worth mentioning that 
$\{\omega_1,\ldots,\omega_m\}$ is the smallest set of characters with the property that 
$\eta_e \eta_{e'} \in \Span(\omega_1,\ldots,\omega_m)$; in other words, no coefficient in the expansion of 
$\eta_e \eta_{e'}$ in terms of~$\omega_1,\ldots,\omega_m$ vanishes. To see this, we revisit the equation
$$\phi_g^*(\eta_e) \psi^*_{\gamma}(\eta_{e'}) = 
\gamma(g^{-1}) \sum_{k=1}^m \xi_k \phi_g^*(\eta_e)(v_k) \omega_k$$
obtained in the second step of the preceding proof. In this equation, the coefficients~$\xi_k$ depend on~$\gamma$, but not on~$g$, whereas the coefficients~$\phi_g^*(\eta_e)(v_k)$ clearly depend on~$g$, but not on~$\gamma$. If one of the coefficients $\xi_1,\ldots,\xi_m$, say~$\xi_i$, would vanish, then this equation would show that for a fixed~$\gamma \in T_e^\perp$, but for all $g \in Q_{e'}^\perp$, we have
$$\phi_g^*(\eta_e) \psi^*_{\gamma}(\eta_{e'}) \in 
\Span(\omega_1,\ldots,\omega_{i-1},\omega_{i+1}\ldots,\omega_m)$$
contradicting the linear independence of the elements on the left. A similar argument with the roles of~$g$ and~$\gamma$ reversed shows that the coefficients~$\phi_g^*(\eta_e)(v_k)$ are all nonzero, so that the coefficient~$\gamma(g^{-1}) \xi_k \phi_g^*(\eta_e)(v_k)$ of~$\omega_k$ in the above equation is nonzero. For $g=1$ and~$\gamma = \eh$, this yields our claim.

We record the following simple consequence of the preceding theorem:
\begin{corollary}
For two primitive idempotents~$e,e' \in \E$, the following conditions are equivalent:
\begin{enumerate}
\item 
$\eta_e \eta_{e'}$ is a character.

\item
$Q_{e'}^\perp \subset T_e$

\item
$\sigma_{A^*,A^*}(\eta_{e} \o \eta_{e'}) = \eta_{e'} \o \eta_{e}$
\end{enumerate}
\end{corollary}
\begin{pf}
In view of the preceding discussion, the first condition is equivalent to the statement that~$m=1$. But since $m=|O|$ by Lemma~\ref{IdealTens}, this is equivalent to the second condition.

It follows directly from the compatibility between product and coproduct discussed in Paragraph~\ref{YetDrinfHopf} that the third condition implies the first. Conversely, it is not difficult to reverse this argument to deduce that the first condition implies the third; a very similar argument can be found in \cite{SoYp}, Prop.~1.1, p.~8. However, it is even easier to use the formula for the quasisymmetry given after Lemma~\ref{FinAbelGroups} to show that the second condition implies the third. From the second condition, we get that also $T_e^\perp \subset Q_{e'}$, and therefore the formula for the quasisymmetry reduces to
\begin{align*}
\sigma_{A^*,A^*}(\eta_{e} \o \eta_{e'}) &= 
\frac{1}{|Q_{e'}^\perp||Q_{e'} \cap T_e^\perp|} \sum_{\gamma \in T_e^\perp} \sum_{g \in Q_{e'}^\perp} \gamma(g^{-1}) \; \psi^*_\gamma(\eta_{e'}) \o \phi^*_g(\eta_{e}) \\
&= \frac{1}{|Q_{e'}^\perp||T_e^\perp|} \sum_{\gamma \in T_e^\perp} \sum_{g \in Q_{e'}^\perp} 
\gamma(g^{-1}) \; \eta_{e'} \o\eta_{e} = \eta_{e'} \o\eta_{e}
\end{align*}
because in this situation we have $\gamma(g^{-1}) = 1$ for~$\gamma \in T_e^\perp$ and~$g \in Q_{e'}^\perp$.
\qed
\end{pf}

\subsection[Ideals and the antipode]{} \label{IdealAntip}
As we will discuss below in greater detail, the image of an ideal under the antipode is in general not again an ideal. However, for a primitive idempotent~\mbox{$e \in \E$}, there is a minimal ideal containing~$e$ whose image under the antipode is again an ideal, namely the ideal~$I_e$: 
\begin{thm}
For a primitive idempotent~$e \in \E$, we have that $\sa(I_e)$ is an ideal. Conversely, if $I \subset A$ is an ideal containing~$e$ with the property that $\sa(I)$ is an ideal, then we have that $I_e \subset I$.
\end{thm}
\begin{pf}
\begin{pflist}
\item
To prove the first assertion, we use the ideal~$\hat{I}_e$ and its complement~$\hat{J}_e$
introduced in Paragraph~\ref{BasicNot}. From the relation between antipode and product stated in Paragraph~\ref{YetDrinfHopf} and the formula for the inverse quasisymmetry in Proposition~\ref{FinAbelGroups}, we get
\begin{align*}
\sa(b) \sa(a) &= \sa(a^\2(\sh^{-1}(a^\1).b)) = \frac{1}{|G|} \sum_{\gamma \in \hat{G}} \sum_{g \in G} 
\gamma(g) \sa(\psi_\gamma(a) \phi_{g}(b))
\end{align*}
Suppose now that $b \in I_e$. If $a \in \hat{J}_e$, then this expression is zero, and in particular in~$\sa(I_e)$. But if $a \in \hat{I}_e$, then we also have that $\psi_\gamma(a) = a$ for all $\gamma \in Q_e$ and therefore can use the formula for the inverse quasisymmetry in Lemma~\ref{FinAbelGroups} instead to get
\begin{align*}
\sa(b) \sa(a) 
= \frac{1}{|Q_e^\perp||Q_e \cap T_e^\perp|} \sum_{\gamma \in T_e^\perp} \sum_{g \in Q_e^\perp} \gamma(g) \sa(\psi_\gamma(a) \phi_{g}(b))
\in \sa(I_e)
\end{align*}
Since $A = \hat{I}_e \oplus \hat{J}_e$, this shows that $\sa(b) \sa(a) \in \sa(I_e)$ for all $a \in A$, which establishes the first assertion because the antipode~$\sa$ is bijective, as we mentioned in Paragraph~\ref{YetDrinfHopf}.

\item
Now suppose that $I \subset A$ is an ideal containing~$e$ with the property that $\sa(I)$ is an ideal. The very last formula yields in the case $b=e$ that 
\begin{align*}
\sa(e) \sa(a) 
= \frac{1}{|Q_e^\perp||Q_e \cap T_e^\perp|} \sum_{\gamma \in T_e^\perp} 
\sum_{g \in Q_e^\perp} \gamma(g) \sa(\psi_\gamma(a) \phi_{g}(e)) 
\end{align*}
for all $a \in \hat{I}_e$. But if $a \in \hat{J}_e$, we saw slightly before that $\sa(e) \sa(a) = 0$. As the right-hand side of the last formula also vanishes in this case, this formula also holds if $a \in \hat{J}_e$, and therefore for all 
$a \in A = \hat{I}_e \oplus \hat{J}_e$. Since $\sa(e) \sa(a) \in \sa(I)$, this implies that
\begin{align*}
\sum_{\gamma \in T_e^\perp} \sum_{g \in Q_e^\perp} \gamma(g) \psi_\gamma(a) \phi_{g}(e) \in I
\end{align*}
for all $a \in A$. Because
$\psi_\gamma(a) \phi_{g}(e) = \phi_{g}(\phi_{g^{-1}}(\psi_\gamma(a)) e) = \eta_e(\phi_{g^{-1}}(\psi_\gamma(a))) \phi_{g}(e)$, this can be written in the form
\begin{align*}
\sum_{g \in Q_e^\perp} \sum_{\gamma \in T_e^\perp}  
\gamma(g) \eta_e(\phi_{g^{-1}}(\psi_\gamma(a))) \phi_{g}(e) \in I
\end{align*}

\item
The elements~$u_h$ introduced in Definition~\ref{IdealTens} are in particular defined when $e=e'$ in the terminology used there, and this is the case which we consider here. For $h \in Q_e^\perp$, we insert $a=\phi_{h^{-1}}(u_h)$ into the condition obtained in the previous step. With the help of Corollary~\ref{IdealTens}, we then find that the element
\begin{align*}
\sum_{g \in Q_e^\perp} &\sum_{\gamma \in T_e^\perp}  
\gamma(g) \eta_e(\phi_{g^{-1}}(\psi_\gamma(\phi_{h^{-1}}(u_h)))) \phi_{g}(e) \\
&= \sum_{g \in Q_e^\perp} \sum_{\gamma \in T_e^\perp}  
\gamma(g) \gamma(h) \eta_e(\phi_{g^{-1}}(\phi_{h^{-1}}(u_h))) \phi_{g}(e) \\
&= |T_e^\perp| \sum_{g \in Q_e^\perp \cap h^{-1} T_e} 
\eta_e(\phi_{g^{-1}h^{-1}}(u_h)) \phi_{g}(e) \\
&= |T_e^\perp| \sum_{g \in h^{-1} Q_e^\perp \cap h^{-1} T_e} \eta_e(u_h) \phi_{g}(e) 
= |T_e^\perp| |Q_e^\perp \cap T_e| \phi_{h^{-1}}(e) 
\end{align*}
is contained in~$I$, and therefore we get that $\phi_{h^{-1}}(e) \in I$ for all 
$h \in Q_e^\perp$, which shows that $I_e \subset I$, as asserted.
\qed
\end{pflist}
\end{pf}

\begin{samepage}
This theorem, which should be compared with~\cite{SoYp}, Lem.~6.3, p.~92, enables us to determine precisely when the image of an ideal under the antipode is again an ideal: 
\begin{corollary}
For an ideal $I \subset A$, the following conditions are equivalent:
\begin{enumerate}
\item 
For all primitive idempotents $e \in I$ and all $g \in Q_e^\perp$, we have $\phi_g(e) \in I$.

\item
$\sa(I)$ is an ideal.
\end{enumerate}
\end{corollary}
\end{samepage}
\begin{pf}
The first condition implies that $I_e \subset I$ for all primitive idempotents $e \in I$, so that
$$I = \sum_{e \in I \cap \E} I_e$$
We then have that
$$\sa(I) = \sum_{e \in I \cap \E} \sa(I_e)$$
is a sum of ideals by the preceding theorem, and therefore itself an ideal.

Conversely, if $\sa(I)$ is an ideal, we have also by the preceding theorem that $I_e \subset I$ for every primitive idempotent $e \in I$, which is equivalent to the first condition.~\qed
\end{pf}

The decisive property of the ideal~$I_e$ is shared by its complement~$J_e$:
\begin{lemma}
For a primitive idempotent $e \in \E$, we have that $\sa(J_e)$ is an ideal.
\end{lemma}
\begin{pf}
We use again the full orbit $\hat{O}_e$, the associated ideal $\hat{I}_e$, and its complement $\hat{J}_e$ from Definition~\ref{BasicNot}. Since~$\hat{J}_e$ is stable under~$G$, the preceding corollary implies that $\sa(\hat{J}_e)$ is an ideal. Now we decompose~$\hat{O}_e$ further as a disjoint union of the form
$$\hat{O}_e = \tilde{O}_1 \cup \cdots \cup \tilde{O}_r$$
where $\tilde{O}_1, \ldots, \tilde{O}_r$ are the orbits of~$Q_e^\perp$ and $\tilde{O}_1=O_e$ is the orbit of~$e$. By construction, all elements of~$\hat{O}_e$ have the same isotropy group as~$e$. If we choose a representative $e_i \in \tilde{O}_i$ for all $i=1,\ldots,r$, where we can choose $e_1=e$, we therefore have
$Q_e^\perp = Q_{e_i}^\perp$, $\tilde{O}_i = O_{e_i}$, and 
$$\hat{I}_e = I_{e_1} \oplus \cdots \oplus I_{e_r}$$
Moreover, $\sa(I_{e_i})$ is an ideal by the preceding theorem. Since we have 
$$J_e = \hat{J}_e \oplus I_{e_2} \oplus \cdots \oplus I_{e_r}$$
we see that $\sa(J_e) = \sa(\hat{J}_e) \oplus \sa(I_{e_2}) \oplus \cdots \oplus \sa(I_{e_r})$ is again an ideal.
\qed
\end{pf}

\subsection[Idempotents and the antipode]{} \label{IdempAntip}
Our considerations so far have not been as symmetric in~$\phi_g$ and~$\psi_\gamma$ as they could have been, as we have used the group~$Q_e$, and not the group~$T_e$, for the definition of the ideal~$I_e$. However, as we pointed out at the end of Paragraph~\ref{FinAbelGroups}, the situation is indeed entirely symmetric in this regard, which in view of Theorem~\ref{IdealAntip} has the following consequence:
\begin{prop}
For a primitive idempotent $e \in \E$, we have
$$\{\phi_g(e) \mid g \in Q_e^\perp\} = \{\psi_\gamma(e) \mid \gamma \in T_e^\perp\}$$
\end{prop}
\begin{pf}
In the same way as in our treatment of Corollary~\ref{IdealTens}, we can consider $A^{\op \cop} = A^{\cop}$ as a left Yetter-Drinfel'd Hopf algebra over~$H^* \cong K[\hat{G}]$. This operation does not change the antipode~$\sa$, so that~$I_e$ is still the minimal ideal~$I$ containing~$e$ with the property that $\sa(I)$ is an ideal. However, this operation exchanges~$\phi_g$ and~$\psi_\gamma$ and consequently~$T_e$ and~$Q_e$, and therefore over~$H^*$ the definition of~$I_e$ would be $I_e := \Span(\{\psi_\gamma(e) \mid \gamma \in T_e^\perp\})$. Since the set of primitive idempotents spanning an ideal is unique, this implies the assertion.
\qed
\end{pf}

It is of course possible to give an alternative proof of this proposition by adapting the proof of Theorem~\ref{IdealAntip} to show that the ideal \mbox{$\Span(\{\psi_\gamma(e) \mid \gamma \in T_e^\perp\})$} is also the minimal ideal~$I$ containing~$e$ with the property that $\sa(I)$ is an ideal. We note that this proposition generalizes \cite{SoYp}, Prop.~6.4, p.~93, where the argument was based on Clifford theory applied to the Radford biproduct. While it is possible to generalize the proof given there to the current situation, the proof that we have given here is considerably simpler.

Obviously, the ideals~$I_e$ and~$\sa(I_e)$ are in general not equal, as these two sets are already usually distinct in the case of ordinary Hopf algebras, i.e., trivial Yetter-Drinfel'd Hopf algebras. However, we can say the following:
\begin{corollary}
$\sa^2(I_e)=I_e$
\end{corollary}
\begin{pf}
In the last step of the proof of Proposition~\ref{Int}, we have already observed that \mbox{$\sa^2 = \theta_A^{-1}$}, the inverse ribbon transformation defined in \cite{SoRib}, Par.~3.5, p.~441. In our situation, the formula given there yields
$$\theta_A^{-1}(a) = a^{\1}.a^{\2}$$
which in view of Corollary~\ref{FinAbelGroups} means for $a \in I_e$ that
$$\theta_A^{-1}(a) = 
\frac{1}{|G|} \sum_{\gamma \in \hat{G}} \sum_{g \in Q_e^\perp} 
\gamma(g^{-1}) \phi_g(\psi_\gamma(a))$$
In the case where $\gamma \notin Q_e T_e^\perp$, we can argue as in the proof of Lemma~\ref{FinAbelGroups} to see that 
$\sum_{g \in Q_e^\perp} \gamma(g^{-1}) \phi_g(\psi_\gamma(a)) = 0$.
The formula above therefore reduces to 
$$\theta_A^{-1}(a) = 
\frac{1}{|G|} \sum_{\gamma \in Q_e T_e^\perp} \sum_{g \in Q_e^\perp} 
\gamma(g^{-1}) \phi_g(\psi_\gamma(a))$$
Now~$I_e$ is stable under~$Q_e^\perp$ by definition and under~$Q_e T_e^\perp$ by the preceding proposition. Therefore, the last equation shows that $\sa^2(a) = \theta_A^{-1}(a) \in I_e$. 
\qed
\end{pf}

The last displayed formula in the preceding proof is not entirely symmetric in~$g$ and~$\gamma$.
However, it would be possible to bring this formula into a symmetric form by using, as in the proof of Lemma~\ref{FinAbelGroups}, the group homomorphism \mbox{$T_e^\perp \times Q_e \rightarrow Q_e T_e^\perp,~(\gamma, \gamma') \mapsto \gamma \gamma'$}. We have not presented this step because it is logically not necessary to prove the preceding corollary.

As the image of an ideal under the antipode is in general not an ideal, we cannot expect that the image of a primitive idempotent under the antipode is again a primitive idempotent. The following lemma describes exactly when this happens:
\begin{lemma}
Suppose that $e \in \E$ is a primitive idempotent. Then $\sa(e)$ is a primitive idempotent if and only if the index of~$e$ is~$1$.
\end{lemma}
\begin{pf}
If $\sa(e)$ is a primitive idempotent, then $I:=Ke$ is an ideal containing~$e$ with the property that $\sa(I)$ is an ideal. Theorem~\ref{IdealAntip} therefore implies that \mbox{$I_e \subset I$}, so that $|G_e| = \dim(I_e) = 1$. Conversely, $|G_e| = 1$ implies $Q_e^\perp = T_e \cap Q_e^\perp$ and therefore $Q_e^\perp \subset T_e$, so that 
$T_e^\perp \subset Q_e$. Now we saw in the second step of the proof of Theorem~\ref{IdealAntip} that
\begin{align*}
\sa(e) \sa(a) 
= \frac{1}{|Q_e^\perp||Q_e \cap T_e^\perp|} \sum_{\gamma \in T_e^\perp} 
\sum_{g \in Q_e^\perp} \gamma(g) \sa(\psi_\gamma(a) \phi_{g}(e)) 
\end{align*}
for all $a \in A$, which then reduces to
\begin{align*}
\sa(e) \sa(a) 
&= \frac{1}{|T_e^\perp|} \sum_{\gamma \in T_e^\perp} 
\sa(\psi_\gamma(a)e) \\
&= \frac{1}{|T_e^\perp|} \sum_{\gamma \in T_e^\perp} 
\eta_e(\psi_\gamma(a)) \sa(e) = \eta_e(a) \sa(e)
\end{align*}
Because the antipode~$\sa$ is bijective, as we noted in Paragraph~\ref{YetDrinfHopf}, this shows on the one hand that $K\sa(e)$ is a one-dimensional ideal, but on the other hand it also shows for $a=e$ that $\sa(e)^2 = \sa(e)$, so that $\sa(e)$ is an idempotent. Since the ideal that it generates is one-dimensional, $\sa(e)$ is primitive (cf.~\cite{Lam}, Prop.~(21.8), p.~320).~\qed
\end{pf}

Our situation has an interesting feature: For the elements~$g \in Q_e^\perp$, the mappings~$\phi_g$ preserve both~$I_e$ and~$\sa(I_e)$ (as do the mappings~$\psi_\gamma$ for~$\gamma \in T_e^\perp$). Inside these ideals, they permute the primitive idempotents, so that we have two permutation representations of~$Q_e^\perp$. Since $\sa \circ \phi_g = \phi_g \circ \sa$ by the $H$-linearity of the antipode, these two representations are isomorphic as linear representations. However, we have seen above that the antipode usually does not map primitive idempotents to primitive idempotents, and is therefore not an isomorphism of permutation representations.

In general, two permutation representations that are isomorphic as linear representations need not be isomorphic as permutation representations. In \cite{Hup1}, Kap.~V, Bsp.~20.10, p.~606, B.~Huppert gives two examples of this phenomenon. In one of these examples, the group is abelian, whereas in the other example, which he attributes to H.~Wielandt, the two permutation representations are transitive, but the group is not abelian. However, in our situation, this phenomenon cannot occur, because two permutation representations of an abelian group that are isomorphic as linear representations are isomorphic as permutation representations if one of them is transitive. This is not very hard to see: If a group element acts as the identity in one representation, then also in the other, because the linear representations are isomorphic. Moreover, this group element has in both permutation representations the same number of fixed points, because this number is just its trace in the corresponding linear representation. But in the transitive permutation representation, a group element that has one fixed point already acts as the identity, because all stabilizers are then conjugate and, as the group is abelian, even equal. We therefore see that this property also holds for the other permutation representation that was not assumed to be transitive, although this one now also turns out to be transitive for cardinality reasons. This shows that both permutation representations are regular in the sense of \cite{Hup1}, Kap.~I, Def.~5.12, p.~27, i.e., isomorphic to the permutation representation on the cosets of their common kernel, and in particular isomorphic to each other (cf.~\cite{Hup1}, Kap.~I, Satz~5.13, p.~27; \cite{KurzwStellm}, Par.~4.1.3, p.~78).

\subsection[The stability set]{} \label{StabSet}
As we have already mentioned in Paragraph~\ref{BasicNot}, we have 
$\check{J}_e^\perp = \Span(\check{U}_e)$. The decisive property of this space is the following:
\begin{thm}
$\Span(\check{U}_e)$ is a Yetter-Drinfel'd Hopf subalgebra of~$A^*$.
\end{thm}
\begin{pf}
We have also pointed out in Paragraph~\ref{BasicNot} that~$\check{U}_e$ is stable under both~$G$ and~$\hat{G}$, so that $\check{J}_e^\perp$ is a submodule and a subcomodule. As it is spanned by group-like elements, $\check{J}_e^\perp$ is in addition a subcoalgebra. It follows from Corollary~\ref{IdealAntip} that $\sa(\check{I}_e)$ and $\sa(\check{J}_e)$ are ideals. Since the primitive idempotents that span $\sa(\check{I}_e)$ are clearly stabilized by~$T_e$ and~$Q_e$, we have $\sa(\check{I}_e) \subset \check{I}_e$ and therefore 
\mbox{$\sa(\check{I}_e) = \check{I}_e$} by dimension considerations. Because
$$A = \sa(\check{I}_e) \oplus \sa(\check{J}_e)= \check{I}_e \oplus \sa(\check{J}_e)$$
this implies that $\sa(\check{J}_e) = \check{J}_e$, which in turn implies that $\sa^*(\check{J}_e^\perp) = \check{J}_e^\perp$.

Since $ \check{U}_e$ evidently contains the counit~$\ea$, it remains to be shown that $\Span(\check{U}_e)$ is multiplicatively closed. To see this, we consider two primitive idempotents~$e'$ and~$e''$ in $\check{O}_e$ and define, as in Paragraph~\ref{IdealTens}, the sets
$$O':= \{\phi_g(e') \mid g \in Q_{e''}^\perp\} \qquad \qquad
O'' := \{\psi_\gamma(e'') \mid \gamma \in T_{e'}^\perp\}$$
Then we know from Lemma~\ref{IdealTens} that $m:=|O'|=|O''|$ and from Theorem~\ref{ProdChar} that there are distinct characters $\omega_1,\ldots,\omega_m$ such that 
$$\Span(\{\eta_{\tilde{e}'} \eta_{\tilde{e}''} \mid \tilde{e}' \in O', \tilde{e}'' \in O'' \}) = 
\Span(\omega_1,\ldots,\omega_m)$$
For $g \in T_e$ and $\tilde{e}' \in O'$, we have $\phi_g(\tilde{e}')=\tilde{e}'$ and therefore
$\phi_g^*(\eta_{\tilde{e}'})=\eta_{\tilde{e}'}$. For exactly the same reason, we have
$\phi_g^*(\eta_{\tilde{e}''})=\eta_{\tilde{e}''}$ for $\tilde{e}'' \in O''$. These two statements together show that $\phi_g^*(\omega_k) = \omega_k$ for all $k=1,\ldots,m$. A very similar argument shows that $\psi_\gamma^*(\omega_k) = \omega_k$ for all $k=1,\ldots,m$ and all $\gamma \in Q_e$, so that $\omega_k \in \check{U}_e$. Since in particular
$$\eta_{e'} \eta_{e''} \in \Span(\omega_1,\ldots,\omega_m) \subset \Span(\check{U}_e)$$
this implies the assertion.
\qed
\end{pf}

By construction, $T_e$ and~$Q_e$ act trivially on~$\check{U}_e$. Therefore, Proposition~\ref{ChangeGroup} implies the following:
\begin{corollary}
$\Span(\check{U}_e)$ is a Yetter-Drinfel'd Hopf algebra over~$K[G_e]$.
\end{corollary}
\begin{pf}
This follows directly from Proposition~\ref{ChangeGroup}, although this proposition was stated for left Yetter-Drinfel'd Hopf algebras, whereas here we have right Yetter-Drinfel'd Hopf algebras. However, we can pass from one to the other with the help of Lemma~\ref{YetDrinfHopf}. 
\qed
\end{pf}

Consequently, primitive idempotents of index~$1$ give rise to trivial Yetter-Drin\-fel'd Hopf subalgebras:
\begin{prop}
For a primitive idempotent~$e \in \E$ of index~$1$, $\check{U}_e$ is a subgroup of the group of units of~$A^*$. Its group ring~$K[\check{U}_e] = \Span(\check{U}_e) \subset A^*$ is therefore a trivial right Yetter-Drinfel'd Hopf subalgebra of~$A^*$. 
\end{prop}
\begin{pf}
If $|G_e|=1$, then the preceding corollary yields that $\Span(\check{U}_e)$ is a Yetter-Drinfel'd Hopf subalgebra over~$K[G_e] \cong K$. The corresponding quasisymmetry is clearly the ordinary flip of tensor factors, so that $\Span(\check{U}_e)$ is trivial, i.e., an ordinary Hopf algebra. It is obviously spanned by the set~$\check{U}_e$ of group-like elements, which therefore form a group.
\qed
\end{pf}

Although the preceding argument is complete, it is possible to give a slightly more explicit proof that $\check{U}_e$ is a group. As we pointed out in the proof of the preceding theorem, $\check{U}_e$ contains the counit~$\ea$. From the proof of Lemma~\ref{IdempAntip}, we know that~$|G_e|=1$ implies $Q_e^\perp \subset T_e$. If~$e' \in \check{O}_e$, so that $T_e \subset T_{e'}$ and $Q_e \subset Q_{e'}$, we therefore have
$$Q_{e'}^\perp \subset Q_e^\perp \subset T_e \subset T_{e'}$$
so that~$e'$ has also index~$1$. Now Lemma~\ref{IdempAntip} yields that $\sa(e')$ is again a primitive idempotent, which is clearly contained in~$\check{O}_e$. This shows that $\sa(\check{O}_e) = \check{O}_e$, and dually we have $\sa^*(\check{U}_e) = \check{U}_e$. Since $\eta_{e'}^{-1}=\sa^*(\eta_{e'})$ for~$e' \in \check{O}_e$, we see that~$\check{U}_e$ contains the inverses of all its elements. 

To show that $\check{U}_e$ is multiplicatively closed, we suppose that~$e'$ and~$e''$ are two primitive idempotents in $\check{O}_e$. We then have
$$Q_{e''}^\perp \subset Q_e^\perp \subset T_e \subset T_{e'}$$
and therefore Corollary~\ref{ProdChar} implies that $\eta_{e'} \eta_{e''}$ is again a character.
Since this character is clearly preserved by~$T_e$ and~$Q_e$, we have 
$\eta_{e'} \eta_{e''} \in \check{U}_e$. From Corollary~\ref{ProdChar}, we also get that
$\sigma_{A^*,A^*}(\eta_{e'} \o \eta_{e''}) = \eta_{e''} \o \eta_{e'}$, which means that $\Span(\check{U}_e)$ is trivial.

\subsection[Dual ideals]{} \label{DualIdeal}
In Paragraph~\ref{ProdChar}, we have investigated the product of two characters belonging to arbitrary primitive idempotents~$e$ and~$e'$. The next step is, exactly as in Paragraph~6.7 of~\cite{SoYp}, to analyze the special case that $e' \in \sa(I_e)$. As we have discussed at the end of Paragraph~\ref{IdempAntip}, we then have $T_e = T_{e'}$, $Q_e = Q_{e'}$, and $I_{e'} = \sa(I_e)$. By Proposition~\ref{IdempAntip}, the orbits~$O$ and~$O'$ introduced at the beginning of Paragraph~\ref{IdealTens} are therefore equal to~$O_e$ and~$O_{e'}$, respectively, and the ideals~$I$ and~$I'$ also introduced there are consequently equal to~$I_e$ and~$I_{e'}$. Proposition~\ref{IdempAntip} also shows that both~$I_e$ and~$I_{e'}$ are invariant under~$T_e^\perp$ and~$Q_e^\perp$.

If $m = |G_e| = |G_{e'}|$ is the index of~$e$ and~$e'$, we have already seen in Theorem~\ref{ProdChar} and the discussion afterwards that there are unique distinct characters $\omega_1,\ldots,\omega_m \in A^*$ with the property that 
$\eta_e \eta_{e'} \in \Span(\omega_1,\ldots,\omega_m)$. By definition, there are corresponding idempotents $e''_1,\ldots,e''_m \in \E$ so that $\omega_k=\eta_{e''_k}$. In our case, when $e' \in \sa(I_e)$, these characters have additional properties, besides those already listed in Theorem~\ref{ProdChar}:
\begin{propb}
\begin{enumerate}
\item
$\{\omega_1,\ldots,\omega_m\} \subset \check{U}_e$

\item 
For $\gamma \in T_e^\perp$, $\{\omega_1,\ldots,\omega_m\}$ is stable under~$\psi_\gamma^*$.

\item
For $g \in Q_e^\perp$, $\{\omega_1,\ldots,\omega_m\}$ is stable under~$\phi_g^*$.

\item
For some $i \le m$, we have $\omega_i = \ea$.

\item
Unless $m=1$, the index of~$\omega_j$ is strictly less than~$m$ for all $j \le m$.
\end{enumerate}
\end{propb}
\begin{pf}
\begin{pflist}
\item
The second and the third assertion follow directly from Theorem~\ref{ProdChar}. This theorem also implies that
$$\Span(U_e U_{e'}) = \Span(\omega_1,\ldots,\omega_m)$$
Since $U_e \subset \check{U}_{e}$ and $U_{e'} \subset \check{U}_{e}$ by construction, Theorem~\ref{StabSet} implies that 
$$\Span(\omega_1,\ldots,\omega_m) \subset \Span(\check{U}_{e})$$
Because distinct group-like elements are linearly independent, as we have already mentioned in the proof of Theorem~\ref{ProdChar}, this yields the first assertion. 

\item
In order to prove the fourth assertion, we recall from Proposition~\ref{BasicNot} that 
\mbox{$\sa^{-1}(e') = \Lma_\1 \eta_{e'}(\Lma_\2)$}, which is contained in~$I_e$ by construction. For the restriction of the module structure on~$I_e$ described in Proposition~\ref{IdealTens} to~$A$ along~$\da$, we therefore have 
\begin{align*}
a.\sa^{-1}(e') &= \da(a).\sa^{-1}(e') = 
a_\1(a_\2{}^\1.\sa^{-1}(e')) \eta_{e'}(a_\2{}^\2) \\
&= a_\1(a_\2{}^\1.\Lma_\1) \eta_{e'}(\Lma_\2) \eta_{e'}(a_\2{}^\2) \\
&= a_\1(a_\2{}^\1.\Lma_\1) \eta_{e'}(a_\2{}^\2 \Lma_\2) \\
&= (a \Lma)_\1 \eta_{e'}((a \Lma)_\2) 
= \ea(a) \Lma_\1 \eta_{e'}(\Lma_\2)
= \ea(a) \sa^{-1}(e')
\end{align*}
From the proof of Theorem~\ref{ProdChar}, we know that the characters 
$\omega_1,\ldots,\omega_m$ arise from the decomposition of this
$A$-module structure on~$I_e$ into one-dimensional submodules. The preceding computation shows that one of these submodules is~$K\sa^{-1}(e')$ and that the corresponding character is~$\ea$. Therefore, we must have~$\omega_i = \ea$ for some $i \le m$, so that $e''_i=\Lma$. 

\item
To prove the fifth assertion, we get from the first assertion that $e''_j \in \check{O}_e$ and hence $T_e \subset T_{e''_j}$ and $Q_e \subset Q_{e''_j}$, which in turn implies that $T_{e''_j}^\perp \subset T_e^\perp$ and $Q_{e''_j}^\perp \subset Q_e^\perp$. From the third assertion, we know that $\{e''_1,\ldots,e''_m\}$ is stable under~$Q_e^\perp$, and therefore
under~$Q_{e''_j}^\perp$. So the orbit~$O_{e''_j}$ of~$e''_j$ under the action of~$Q_{e''_j}^\perp$ is contained in $\{e''_1,\ldots,e''_m\}$, which shows that the index~$|O_{e''_j}|$ of~$e''_j$ is less than or equal to~$m$. To show that it is strictly less than~$m$, we distinguish two cases: If $j=i$, we have $e''_i=\Lma$, which has index~$1$, and~$1$ is less than~$m$ by assumption. If $j \neq i$, then $e''_i \notin O_{e''_j}$, 
so~$|O_{e''_j}| < m$, as asserted.
\qed
\end{pflist}
\end{pf}

This proposition shows that after renumbering the characters, we can, and will, assume that 
$\omega_1 = \ea$, so that $e''_1 = \Lma$. Its proof shows that for the corresponding eigenvector~$v_1$ introduced in Paragraph~\ref{ProdChar}, we can choose $v_1=\sa^{-1}(e')$.

As pointed out at the beginning of this paragraph, we have \mbox{$\sa(I_e) = I_{e'}$}. The statement $\sa^2(I_e) = I_e$ in Corollary~\ref{IdempAntip} therefore means that 
$\sa(I_{e'}) = I_e$. In addition, we get from $A = I_e \oplus J_e$ that 
$$A = \sa(I_e) \oplus \sa(J_e) = I_{e'} \oplus \sa(J_e)$$
Since $\sa(J_e)$ is an ideal by Lemma~\ref{IdealAntip}, it must be the unique ideal that complements~$I_{e'}$; in other words, we have $\sa(J_e) = J_{e'}$. A similar argument shows that 
$\sa(J_{e'}) = J_e$. These observations are used in the proof of the following lemma:
\begin{lemma}
$\Span(\omega_1,\ldots,\omega_m)$ is stable under the antipode~$\sa^*$.
\end{lemma}
\begin{pf}
We have already noted in Paragraph~\ref{BasicNot} that 
$$J_e^\perp = \Span(U_e) \qquad \qquad J_{e'}^\perp = \Span(U_{e'})$$
and the discussion above implies that $\sa^*(J_{e'}^\perp)=J_e^\perp$ and $\sa^*(J_e^\perp)=J_{e'}^\perp$. Moreover, we know from the proof of the preceding proposition that
$$\Span(U_e U_{e'}) = \Span(\omega_1,\ldots,\omega_m)$$
Let us now suppose that $\tilde{e} \in O_e$ and $\tilde{e}' \in O_{e'}$. By Proposition~\ref{YetDrinfHopf}, $A^*$ is a right Yetter-Drinfel'd Hopf algebra over~$H$ with antipode~$\sa^*$, and from the relation between antipode and multiplication also stated in Paragraph~\ref{YetDrinfHopf} together with the formulas for the quasisymmetry given after Lemma~\ref{FinAbelGroups} we get
$$\sa^*(\eta_{\tilde{e}} \eta_{\tilde{e}'}) = 
\frac{1}{|Q_e^\perp||Q_e \cap T_e^\perp|} 
\sum_{\gamma \in T_e^\perp} \sum_{g \in Q_e^\perp}
\gamma(g^{-1}) \sa^*(\psi_{\gamma}^*(\eta_{\tilde{e}'})) \sa^*(\phi^*_{g}(\eta_{\tilde{e}}))$$
For $\gamma \in T_e^\perp$ and $g \in Q_e^\perp$, we have
$\psi_{\gamma}^*(\eta_{\tilde{e}'}) \in J_{e'}^\perp$ and 
$\phi^*_{g}(\eta_{\tilde{e}}) \in J_e^\perp$, because $J_e$ and~$J_{e'}$ are stable under $Q_e^\perp = Q_{e'}^\perp$ by construction and under $T_e^\perp = T_{e'}^\perp$ by Proposition~\ref{IdempAntip}. Therefore, the above formula shows that 
$\sa^*(\eta_{\tilde{e}} \eta_{\tilde{e}'})$ is contained in the span of 
$J_e^\perp J_{e'}^\perp$, which is $\Span(\omega_1,\ldots,\omega_m)$, as asserted.
\qed
\end{pf}

We note that, since $\sa^*(\omega_k) = \omega_k^{-1}$, this lemma can also be stated in the form
$$\Span(\omega_1^{-1},\ldots,\omega_m^{-1})=\Span(\omega_1,\ldots,\omega_m)$$

The spaces $J_e^\perp$ and~$J_{e'}^\perp$ just considered are preserved under left and right multiplication by~$\omega_k$:
\begin{thm}
For all $k=1,\ldots,m$, we have
$$\omega_k J_{e}^\perp = J_{e}^\perp \omega_k = J_{e}^\perp \qquad  \qquad
\omega_k J_{e'}^\perp = J_{e'}^\perp \omega_k = J_{e'}^\perp$$
Moreover, we have $\Span(U_{e'} U_e) = \Span(\omega_1,\ldots,\omega_m)$.
\end{thm}
\begin{pf}
\begin{pflist}
\item
For $\tilde{e} \in O_e$, we have that
$\eta_{\tilde{e}}^{-1} = \sa^*(\eta_{\tilde{e}}) \in \sa^*(J_{e}^\perp) = J_{e'}^\perp$ and therefore
$$J_{e}^\perp \eta_{\tilde{e}}^{-1} \subset J_{e}^\perp J_{e'}^\perp
\subset \Span(\omega_1,\ldots,\omega_m)$$
Since the spaces on the left and on the right have the same dimension, we get
$J_{e}^\perp \eta_{\tilde{e}}^{-1} =  \Span(\omega_1,\ldots,\omega_m)$ and consequently
$J_{e}^\perp = \Span(\omega_1 \eta_{\tilde{e}},\ldots,\omega_m \eta_{\tilde{e}})$. This shows
that $\omega_k \eta_{\tilde{e}} \in J_{e}^\perp$, which yields $\omega_k J_{e}^\perp = J_{e}^\perp$.

\item
Similarly, we have for $\tilde{e}' \in O_{e'}$ that 
$\eta_{\tilde{e}'}^{-1} = \sa^*(\eta_{\tilde{e}'}) \in \sa^*(J_{e'}^\perp) = J_{e}^\perp$ and therefore
$$\eta_{\tilde{e}'}^{-1} J_{e'}^\perp \subset J_{e}^\perp J_{e'}^\perp
\subset \Span(\omega_1,\ldots,\omega_m)$$
Again by comparing dimensions, we get
$\eta_{\tilde{e}'}^{-1} J_{e'}^\perp = \Span(\omega_1,\ldots,\omega_m)$ and consequently
$J_{e'}^\perp = \Span(\eta_{\tilde{e}'} \omega_1,\ldots,\eta_{\tilde{e}'} \omega_m)$. This shows
that $\eta_{\tilde{e}'} \omega_k \in J_{e'}^\perp$, which yields 
\mbox{$J_{e'}^\perp \omega_k = J_{e'}^\perp$}.

\item 
By applying Theorem~\ref{ProdChar} in the case where the idempotents are equal, we see that 
$\Span(U_e U_e)$ is an $m$-dimensional space. Hence we have
$$J_{e}^\perp \eta_e = \Span(U_e \eta_e) = \Span(\eta_e U_e) = \eta_e J_{e}^\perp$$
Since we saw in the first step that $\Span(\omega_1 \eta_e,\ldots,\omega_m \eta_e) = J_{e}^\perp$,
it follows that
$\Span(\eta_e \omega_1 \eta_e,\ldots, \eta_e \omega_m \eta_e) = \eta_e J_{e}^\perp = J_{e}^\perp \eta_e$.

\item
As we saw above, we have not only $I_{e'} = \sa(I_e)$, but also $I_e = \sa(I_{e'})$. We can therefore interchange the roles of~$e$ and~$e'$ and find characters~$\omega'_1,\ldots,\omega'_m$ with the property that
$$\Span(U_{e'} U_e) = \Span(\omega'_1,\ldots,\omega'_m)$$
which by the first two steps satisfy 
$\omega'_k J_{e'}^\perp = J_{e'}^\perp$ and $J_{e}^\perp \omega'_k = J_{e}^\perp$, and by the second step also $\Span(\eta_e \omega'_1,\ldots, \eta_e \omega'_m)= J_{e}^\perp$. This implies
\begin{align*}
\Span(\eta_e \omega'_1 \eta_e,\ldots, \eta_e \omega'_m \eta_e) = J_{e}^\perp \eta_e
\end{align*}

\item
Comparing the last two steps, we get 
$$\Span(\eta_e \omega_1 \eta_e,\ldots,\eta_e \omega_m \eta_e) = 
\Span(\eta_e \omega'_1 \eta_e,\ldots, \eta_e \omega'_m \eta_e)$$
which in turn clearly implies $\Span(\omega_1,\ldots,\omega_m) = \Span(\omega'_1,\ldots, \omega'_m)$ by multiplying with~$\eta_e^{-1}$ on the left and on the right. As we have already recalled in the proof of the proposition above, distinct group-like elements are linearly independent, and hence we must have 
$\{\omega_1,\ldots,\omega_m\} = \{\omega'_1,\ldots, \omega'_m\}$. For a given $k=1,\ldots,m$, we can therefore find $l=1,\ldots,m$ so that $\omega_k = \omega'_l$, and then have
$\omega_k J_{e'}^\perp = \omega'_l J_{e'}^\perp = J_{e'}^\perp$ and 
$J_{e}^\perp \omega_k = J_{e}^\perp \omega'_l =J_{e}^\perp$, as asserted.
\qed
\end{pflist}
\end{pf}

This theorem obviously implies that multiplication by the inverse of~$\omega_k$ also preserves these spaces; i.e., we have
$$\omega_k^{-1} J_{e}^\perp = J_{e}^\perp \omega_k^{-1} = J_{e}^\perp \qquad 
\qquad
\omega_k^{-1} J_{e'}^\perp = J_{e'}^\perp \omega_k^{-1} = J_{e'}^\perp$$

\subsection[The core]{} \label{Core}
In the last paragraph, we have considered two primitive idempotents~$e$ and~$e'$ that satisfy the restriction $e' \in \sa(I_e)$, and we have obtained distinct characters $\omega_1,\ldots,\omega_m$ with the property that
$$\Span(U_e U_{e'}) = \Span(U_{e'} U_e) = \Span(\omega_1,\ldots,\omega_m)$$
where $m = |G_e| = |G_{e'}|$. For a given primitive idempotent $e \in \E$, this space does not depend on the choice of the primitive idempotent~$e' \in \sa(I_e)$, because two different choices of such idempotents nonetheless generate the same orbit~$O_{e'}$. Therefore, the following definition is meaningful:
\begin{defn}
For a primitive idempotent~$e \in \E$, the space~$\Span(U_e U_{e'})$, where~$e'$ is some primitive idempotent in~$\sa(I_e)$, is called the core of~$e$, or alternatively the core of~$\eta_e$.
\end{defn}

As recalled in the preceding paragraph, the dual space $A^*$ is a right Yetter-Drinfel'd Hopf algebra over~$H$. As we will see in a moment that the core $\Span(\omega_1,\ldots,\omega_m)$ is multiplicatively closed, it seems reasonable to conjecture that the core is a Yetter-Drinfel'd Hopf subalgebra of~$A^*$. This, however, cannot be the case, as it is not invariant under the entire group~$G$, but only, as we have seen in Proposition~\ref{DualIdeal}, under the subgroup~$Q_e^\perp$. Nevertheless, it is a Yetter-Drinfel'd Hopf subalgebra over a different group:
\begin{thm}
The core $\Span(\omega_1,\ldots,\omega_m)$ is a Yetter-Drinfel'd Hopf subalgebra of the right $K[G_e]$-Yetter-Drinfel'd Hopf algebra~$\Span(\check{U}_e)$.
\end{thm}
\begin{pf}
\begin{pflist}
\item
First, we have to verify that the core is multiplicatively closed. But this follows from Theorem~\ref{DualIdeal}, as we have
\begin{align*}
(U_e U_{e'}) (U_e U_{e'}) &\subset U_e \Span(\omega_1,\ldots,\omega_m) U_{e'} \\
&\subset J_{e}^\perp \Span(\omega_1,\ldots,\omega_m) J_{e'}^\perp
\subset J_{e}^\perp J_{e'}^\perp \subset \Span(\omega_1,\ldots,\omega_m)
\end{align*}

\item
From the first assertion in Proposition~\ref{DualIdeal}, we know that the core is contained in $\Span(\check{U}_e)$, and from its fourth assertion we know that the core contains the counit~$\ea$ and is therefore a unital subalgebra of~$\Span(\check{U}_e)$. Because $\omega_1,\ldots,\omega_m$ are group-like elements of the dual, the core is also a subcoalgebra. Again by Proposition~\ref{DualIdeal}, our space is stable under the maps~$\phi_g^*$ and $\psi_\gamma^*$ for $g \in Q_e^\perp$ and $\gamma \in T_e^\perp$. The action of~$Q_e^\perp$ clearly factors over \mbox{$G_e = Q_e^\perp/(T_e \cap Q_e^\perp)$}, which means that the core is a right $K[G_e]$-submodule of~$\Span(\check{U}_e)$. Similarly, the action of~$T_e^\perp$ factors over~$T_e^\perp/(Q_e \cap T_e^\perp)$, and we saw in Lemma~\ref{ChangeGroup} that this group is isomorphic to the character group~$\hat{G}_e$ of~$G_e$. From the way how the coaction was constructed in Paragraph~\ref{ChangeGroup}, we see that the fact that the core is a $K[T_e^\perp/(Q_e \cap T_e^\perp)]$-submodule means that it is also a $K[G_e]$-subcomodule. Finally, Lemma~\ref{DualIdeal} implies that $\Span(\omega_1,\ldots,\omega_m)$ is stable under the antipode and hence a Yetter-Drinfel'd Hopf subalgebra of~$\Span(\check{U}_e)$.
\qed
\end{pflist}
\end{pf}

As we have pointed out above, the core is not a Yetter-Drinfel'd Hopf subalgebra of~$A^*$. However, M.~Takeuchi has introduced a version of the notion of a braided Hopf algebra which does not use the quasisymmetry of some outside category, but rather lists this map as part of the data of the braided Hopf algebra (cf.~\cite{TakSurvBraid}, Def.~5.1, p.~310). Since they are both Yetter-Drinfel'd Hopf algebras, although over different Hopf algebras, both the core and~$A^*$ are Hopf algebras in Takeuchi's sense (cf.~\cite{TakSurvBraid}, Thm.~5.7, p.~314). Furthermore, the proof of Proposition~\ref{ChangeGroup} shows that the core is a Yang-Baxter subspace in Takeuchi's terminology (cf.~\cite{TakSurvBraid}, Eq.~(6.3), p.~314). However, in this terminology, it is not a categorical Yang-Baxter subspace (cf.~\cite{TakSurvBraid}, Eq.~(6.4), p.~315).

As Takeuchi explains, Scharfschwerdt's version of the Nichols-Zoeller theorem cited in Paragraph~\ref{Int} can be used to establish freeness for categorical braided Hopf subalgebras in the sense just discussed (cf.~\cite{TakSurvBraid}, Thm.~7.3, p.~316). As the core is not categorical, this result does not apply directly to the core viewed as a subalgebra of~$A^*$. Nonetheless, there are two ways to arrive at the same conclusion:
\begin{corollary}
$A^*$ is free as a left and right module over $\Span(\omega_1,\ldots,\omega_m)$. In particular, the index~$m$ of~$e$ divides~$\dim(A)$.
\end{corollary}
\begin{pf}
We have already seen in the proof of the preceding theorem that the core is a unital subalgebra, and also a subcoalgebra, of~$A^*$. Although the module and the comodule structure of both the core and~$A^*$ are on the right, Lemma~\ref{YetDrinfHopf} implies that Theorem~\ref{Int} also holds for right Yetter-Drinfel'd Hopf algebras, and so this theorem immediately yields the assertion.
\qed
\end{pf}

\begin{sloppypar}
However, there is a second proof of this corollary that does not rely on the refined version of the Nichols-Zoeller theorem given in Theorem~\ref{Int}, but only requires the version of this theorem for Yetter-Drinfel'd Hopf algebras by B.~Scharf\-schwerdt mentioned above and cited in Paragraph~\ref{Int}: By Theorem~\ref{StabSet}, $\Span(\check{U}_e)$ is a Yetter-Drinfel'd Hopf subalgebra of~$A^*$, both considered as Yetter-Drinfel'd Hopf algebras over~$K[G]$. Therefore, $A^*$ is free as a left and right module over $\Span(\check{U}_e)$. But by Corollary~\ref{StabSet}, $\Span(\check{U}_e)$ can also be considered as a Yetter-Drinfel'd Hopf algebra over~$K[G_e]$, and by the preceding theorem the core is a Yetter-Drinfel'd Hopf subalgebra of $\Span(\check{U}_e)$, both considered as Yetter-Drinfel'd Hopf algebras over~$K[G_e]$. Therefore, $\Span(\check{U}_e)$ is free as a left and right module over $\Span(\omega_1,\ldots,\omega_m)$. Since freeness is transitive, this proves the preceding corollary.
\end{sloppypar}

\subsection[The triviality theorems]{} \label{TrivThm}
We are now prepared to derive our main results. Before we state them, we first recall that we are considering a commutative semisimple left Yetter-Drinfel'd Hopf algebra over the group ring~$K[G]$ of a finite abelian group~$G$. The base field~$K$ is assumed to be algebraically closed of characteristic zero. In this situation, we can say the following:
\begin{thm}
If~$\dim(A)$ and~$|G|$ are relatively prime, then $A$ is trivial.
\end{thm}
\begin{pf}
\begin{pflist}
\item
The index $|G_e|$ of a primitive idempotent~$e$ divides~$|G|$, and we have shown in Corollary~\ref{Core} that $|G_e|$ divides~$\dim(A)$. Hence our assumption implies that every primitive idempotent has index~$1$. Proposition~\ref{StabSet} therefore yields that~$\check{U}_e$ is a group and that its group ring~$K[\check{U}_e] = \Span(\check{U}_e)$ is a trivial Yetter-Drinfel'd Hopf subalgebra of~$A^*$.

\item
We now consider two primitive idempotents~$e$ and~$e'$ of~$A$ that are not necessarily related in any way. It follows from Theorem~\ref{ProdChar} that
$$\Span(\{\phi_g^*(\eta_{e}) \eta_{e'} \mid g \in Q_{e'}^\perp\}) = 
\Span(\{\eta_{e} \psi_\gamma^*(\eta_{e'})  \mid \gamma \in T_{e}^\perp\})$$
because the spanning sets on both sides of this equation are linearly independent and have by Lemma~\ref{IdealTens} the same number of elements. Multiplying by~$\eta_{e}^{-1}$ from the left and by~$\eta_{e'}^{-1}$ from the right, we get that
$$\Span(\{\eta_{e}^{-1} \phi_g^*(\eta_{e}) \mid g \in Q_{e'}^\perp\}) = 
\Span(\{\psi_\gamma^*(\eta_{e'}) \eta_{e'}^{-1} \mid \gamma \in T_{e}^\perp\})$$
As the left-hand side is a subset of~$K[\check{U}_e]$ and the right-hand side is a subset of~$K[\check{U}_{e'}]$, the spanning sets on both sides are now group-like elements. Since distinct group-like elements in any coalgebra are linearly independent, as we already mentioned in the proof of Theorem~\ref{ProdChar}, we get the stronger statement that
$$\{\eta_{e}^{-1} \phi_g^*(\eta_{e}) \mid g \in Q_{e'}^\perp\} = 
\{\psi_\gamma^*(\eta_{e'}) \eta_{e'}^{-1} \mid \gamma \in T_{e}^\perp\}$$
In particular, for every $g \in Q_{e'}^\perp$ there exists $\gamma \in T_{e}^\perp$ such that
$$\eta_{e}^{-1} \phi_g^*(\eta_{e}) = \psi_\gamma^*(\eta_{e'}) \eta_{e'}^{-1}$$
If we introduce the element $\eta':= \psi_\gamma^*(\eta_{e'}) \eta_{e'}^{-1} \in \check{U}_{e'}$, which is invariant 
under~\mbox{$g \in Q_{e'}^\perp \subset T_{e'}$}, we can write the last equation in the form 
$\phi_g^*(\eta_{e}) = \eta_{e} \eta'$. More generally, we have for all nonnegative integers~$i \in \N_0$ that
$$\phi_{g^i}^{*}(\eta_{e}) = \eta_{e} \eta'^i$$
This is obvious for~$i=0$, holds for $i=1$ by construction, and then follows inductively, since
$$\phi_{g^{i+1}}^{*}(\eta_{e}) = \phi_{g}^{*}(\eta_{e} \eta'^i) = \phi_{g}^{*}(\eta_{e}) \phi_{g}^{*}(\eta'^i)
= \phi_{g}^{*}(\eta_{e}) \eta'^i = \eta_{e} \eta'^{i+1}$$

\item
Now let $k$ be the order of~$g \in G$, so that $g^k=1$. Then we have
$$\eta_{e} = \phi_{g^k}^{*}(\eta_{e}) = \eta_{e} \eta'^k$$ 
and therefore $\eta'^k=1$. So the order of~$\eta'$ divides~$k$, which divides the order of~$G$. But the order of~$\eta'$ also divides the order of~$\check{U}_{e'}$, which is the dimension of~$K[\check{U}_{e'}]$. The dimension of~$K[\check{U}_{e'}]$ in turn divides the dimension of~$A$ by the Nichols-Zoeller theorem for Yetter-Drinfel'd Hopf algebras cited in Paragraph~\ref{Int} (the refined version given in Theorem~\ref{Int} is not necessary here, although it also yields the result). Because by assumption these numbers are relatively prime, the order of~$\eta'$ must be~$1$. Therefore we have \mbox{$\eta'=1$} and $\phi_g^*(\eta_{e}) = \eta_{e}$. This shows that $\eta_{e}$ is invariant under~$Q_{e'}^\perp$, which means that 
$Q_{e'}^\perp \subset T_e$. Now Corollary~\ref{ProdChar} implies that 
$$\sigma_{A^*,A^*}(\eta_{e} \o \eta_{e'}) = \eta_{e'} \o \eta_{e}$$
Since the characters form a basis of~$A^*$, this shows that~$A^*$ is trivial as a right Yetter-Drinfel'd Hopf algebra. But as we know from Paragraph~\ref{RightYetMod} that~$\sigma_{A,A}$ and~$\sigma_{A^*,A^*}$ are adjoints of each other, this implies that~$A$ is trivial, as asserted.~\qed
\end{pflist}
\end{pf}

As already mentioned in the introduction, this result generalizes \cite{SoYp}, Cor.~6.7, p.~100. It may be worth pointing out that the assumption that~$K$ is algebraically closed is in fact not necessary: According to Definition~\ref{YetDrinfHopf}, triviality means that the equation \mbox{$\sigma_{A,A}(a \o a') = a' \o a$} holds for all $a, a' \in A$. The validity of this equation does not depend on the base field. We can therefore enlarge the base field to its algebraic closure to decide this question, because our assumptions will still be satisfied over the enlarged base field. This holds in particular for the semisimplicity assumption, as we discussed in Paragraph~\ref{Int}.

We now return to the situation where~$K$ is an algebraically closed field of characteristic zero, but consider instead a finite-dimensional cocommutative cosemisimple right Yetter-Drinfel'd Hopf algebra~$A$ over the group ring~$K[G]$ of a finite abelian group~$G$. Although~$A$ itself need not be trivial, it contains at least a trivial part:
\begin{prop}
If $\dim(A)>1$, then $A$ contains a trivial Yetter-Drinfel'd Hopf subalgebra~$B$ with $\dim(B)>1$.
\end{prop}
\begin{pf}
If we can find a primitive idempotent~$e \in A^*$ of index~$1$ that is different from the integral, then Proposition~\ref{StabSet} yields that $B:= K[\check{U}_e] \subset A^{**} \cong A$ is a trivial Yetter-Drinfel'd Hopf subalgebra. As $B$ contains both the counit and~$\eta_e$, we have $\dim(B)>1$.

In the case where every primitive idempotent has index~$1$, we can obviously find such an idempotent. Otherwise there are primitive idempotents whose index is different from~$1$, and among those we choose a primitive idempotent~$e'$ of minimal index $m:=|G_{e'}|>1$. If then~$e$ is a primitive idempotent that corresponds to a character in the core of~$e'$, we have $|G_{e}|<m$ by Proposition~\ref{DualIdeal}, and therefore $|G_{e}|=1$. Since there are~$m$ such idempotents, one of them is different from the integral, which in view of the discussion above establishes our assertion.
\qed
\end{pf}

It should be noted that the Nichols-Zoeller theorem for Yetter-Drinfel'd Hopf algebras cited in Paragraph~\ref{Int} implies in this situation that the dimension of~$B$ divides the dimension of~$A$.

\end{document}